\documentclass[a4paper,11pt]{article}
\usepackage{amssymb}
\usepackage{color}
\usepackage{graphicx}
\usepackage{amsmath}

\newcommand{\R}{\mathbb{R}}

\newcommand{\N}{\mathbb{N}}
\newcommand{\cL}{\mathcal{L}}

\newcommand{\cuad}{{\sqcap\kern-.68em\sqcup}}

\newcommand{\s}{s}

 \newcommand{\cI}{\mathcal{I}}
 
\newcommand{\cC}{\mathcal{C}}

\newcommand{\cD}{\mathcal{D}}

\newcommand{\cO}{\mathcal{O}}

\newcommand{\cF}{\mathcal{F}}

\newtheorem{theorem}{Theorem}[section]
\newtheorem{proposition}{Proposition}[section]

\newtheorem{lemma}{Lemma}[section]
\newtheorem{corollary}{Corollary}[section]
\newtheorem{remark}{Remark}[section]
\newcommand{\bremark}{\begin{remark} \em}
\newcommand{\eremark}{\end{remark} }

\newcommand{\cP}{\mathcal{P}}
\newcommand{\bP}{\mathbb{P}}

\newcommand{\ms}{\medskip}

\newcommand{\loc}{{\rm loc}}

 \allowdisplaybreaks

\usepackage{cite}
\begin{document}

\begin{center}{\bf  \large   On  positive solutions of  critical semilinear  equations\\[2mm]
 involving the Logarithmic Laplacian

 } \bigskip\bigskip

{\small

  Huyuan Chen\footnote{chenhuyuan@yeah.net}

\bigskip
Shanghai Institute for Mathematics and Interdisciplinary Sciences, \\
Fudan University, \  Shanghai ‌200433, PR  China \\[6pt]

The University of Sydney, School of Mathematics and Statistics,
\\ NSW 2006,
Australia \\[16pt]

Feng Zhou\footnote{fzhou@math.ecnu.edu.cn}\bigskip

 CPDE, School of Mathematical Sciences, East China Normal University, \\
 Shanghai 200241, PR China\\[6pt]
  NYU-ECNU Institute of Mathematical Sciences at NYU-Shanghai,\\
Shanghai 200120, PR China}  \\[18pt]

\begin{quote}

{\bf Abstract.} In this paper, we  classify the   solutions of the critical semilinear problem  involving the logarithmic Laplacian
$$(E)\qquad \qquad\qquad\qquad\qquad \cL_\Delta u= k u\ln u,\qquad u\geq0 \quad \  {\rm in}\ \   \R^n, \qquad\qquad\qquad\qquad\qquad\qquad$$
where   $k\in(0,+\infty)$,  $\cL_\Delta$  is the logarithmic Laplacian in $\R^n$ with $n\in\N$, and
$s\ln s=0$ if $s=0$.
When $k=\frac4n$, problem $(E)$ only has the solutions with the form
$$u_{\tilde x,t}(x)=\beta_n \Big(\frac{t}{t^2+|x-\tilde x|^2)}\Big)^{\frac{n}{2}}\quad \text{  for any $t>0$,  $\tilde x\in\R^n$},$$ where   $\beta_n=2^{\frac n2} e^{\frac n2\psi(\frac n2) }>0$.  When $k\in(0,+\infty)\setminus\{\frac 4n\}$,
 problem $(E)$  has no any positive solution under some suitable assumptions.
\medskip


\end{quote}

  \end{center}



\setcounter{equation}{0}
\section{Introduction}
Let $\cL_\Delta$ be  the logarthmic Laplacian  in $\R^n$ with $n\in\N$, the set of positive integers, initiated in \cite{CW0},   defined by
$$
\cL_\Delta  u(x)= c_{n} \int_{\R^n  } \frac{ u(x)1_{B_1(x)}(y)-u(y)}{|x-y|^{n} } dy + \rho_n u(x),
$$
where
$$
c_n\!:= \pi^{-  n/2}  \Gamma( n/2) = \frac{2}{\omega_{{n-1}}}, \qquad \rho_n\!:=2 \ln 2 + \psi(\tfrac{n}{2}) -\gamma,
$$
$\omega_{{n-1}}:=|\mathbb{S}^{n-1}|$, the area of the unit sphere of $\R^{n}$, and $\gamma= -\Gamma'(1)$ is the Euler Mascheroni constant, $\psi = \frac{\Gamma'}{\Gamma}$ is the digamma function associated to the Gamma function $\Gamma$.
Our aim of this article is to classify the positive  solutions for the critical semilinear problem  involving the logarithmic Laplacian
\begin{equation}\label{eq 1.1}
\left\{ \arraycolsep=1pt
\begin{array}{lll}
\cL_\Delta u= \frac4n u\ln u \quad \  &{\rm in}\ \   \R^n,\\[2mm]
 \phantom{ \cL_{\Delta}   }
  u\gneqq 0\quad \ &{\rm{in}}  \ \ \R^n,
\end{array}
\right.
\end{equation}
where   the function  $t\to t\ln t $ is defined on $[0,+\infty)$ with  $t\ln t=0$ for  $t=0$.

In the last decades, the seminlinear  equations  involving integro-differential operators have been studied extensively and widely
 inspired by a renewed and increasing interest and by their wide range of applications. Also important progresses have been made in
understanding nonlocal phenomena from a PDE or a probabilistic point of view, see e.g.
\cite{CR,CS0,CS1,CF,RS1,RS,Si,musina-nazarov,L,T}.
Among nonlocal differential order operators the simplest examples are the fractional powers of the Laplacian, which exhibit  many phenomenological  properties.
Recall that for $\s\in(0,1)$, the fractional Laplacian
$(-\Delta)^\s$ can be written as a singular integral operator defined in the principle value sense
$$
 (-\Delta)^s  u(x)=c_{n,s} \lim_{\epsilon\to0^+} \int_{\R^n\setminus B_\epsilon(x) }\frac{u(x)-
u(y)}{|x-y|^{n+2s}}  d y ,
$$
where $c_{n,s}=2^{2\s}\pi^{-\frac n2}\s\frac{\Gamma(\frac{n+2\s}2)}{\Gamma(1-\s)}>0$ is a normalized constant such that  for a function $u \in C^\infty_c(\R^n)$,
\begin{equation*}
  \label{eq:Fourier-representation}
\cF((-\Delta)^s u)(\xi):=(2\pi)^{-\frac{n}{2}}\int_{\R^n}e^{-{\rm i}x\cdot \xi}((-\Delta)^s u)(x) d x= |\xi|^{2s}\widehat u (\xi)\quad \text{for all $\xi \in \R^n$}.
\end{equation*}
Here and in the sequel both $\mathcal{F}$ and $\widehat \cdot$ denote the Fourier transform.   Two  known limits hold when $\s$ tends $0$ or $1$:
\begin{equation*}
\lim_{\s\to1^-}(-\Delta)^\s u(x)=-\Delta u(x)
\quad \text{and}\quad \lim_{\s\to0^+}  (-\Delta)^\s  u(x) = u(x)\quad\ \text{for $u\in \mathcal{C}^2_c(\R^n)$,}
\end{equation*}
see e.g. \cite{DPV}.
Furthermore  the following surprising expansion at $s=0$ is  proved in \cite{CW0}
$$
(-\Delta)^s  u(x) = u(x) + s \cL_\Delta u (x) + o(s) \quad \text{as }\,s\to 0^+\quad\text{ for all  }u \in \cC^2_c(\R^n) \text{ and }
x\in \R^n, $$
 where the formal operator
$$
\cL_\Delta:= \frac{d}{d\s}\Big|_{\s=0} (-\Delta)^\s
$$
 is the {\em logarithmic Laplacian}; more precisely (see e.g. \cite{CW0}),
\begin{enumerate}
\item[(i)] for $1 < p \le \infty$, there holds $\cL_\Delta  u \in L^p(\R^n)$ and $\frac{(-\Delta)^\s u- u}{\s} \to \cL_\Delta  u$ in $L^p(\R^n)$ as $\s \to 0^+$;
\item[(ii)] $\mathcal{F}(\cL_\Delta u)(\xi) = (2 \ln |\xi|)\,\widehat u (\xi)$ 
 \, for a.e. $\xi \in \R^n$.
\end{enumerate}
Recently, the study  involving the logarithmic Laplacian has been investigated in various models: the eigenvalues estimates \cite{CV,LW},  Poisson problem \cite{ClS,EV,HS,HRS,JSW}, Heat equation \cite{CDK,CV1},
application to the analysis of the Levy Flight \cite{DGV} and the references therein. \smallskip

  The natural domain of definition of $\cL_\Delta$  is the set of uniformly Dini continuous functions  $u$ in $\R^n$
with the integral restriction
$$\|u\|_{L^1_0(\R^n)}:=\int_{\R^n}\frac{|u(x)|}{1+|x|^{n}}dx<\infty.
$$
Precisely,  by denoting the module of continuity of $u$ at $x\in\R^n$ by
$$
\varpi_{u,x }: [0,1] \to [0,\infty),\qquad  \varpi_{u,x}(r)= \sup_{ y \in \R^n,\, |y-x|\le r} |u(y)-u(x)|,
$$
 a function $u$ is said to be  Dini continuous  at $x$ if
$$\displaystyle \int_0^1 \frac{\varpi_{u,x}(r)}{r}  dr < +\infty.$$
Let  $\cD^0(\R^n)$ be the set of all functions  $u$ having the Dini continuity at any $x\in\R^n$. {\it Here a nonnegative function $u$ is said to be a solution of \eqref{eq 1.1} if $u\in L^1_0(\R^n)\cap \cD^0(\R^n)$ satisfies the equation pointwisely.  }

 \smallskip

Our first result is on the existence and classification of positive solutions of \eqref{eq 1.1}.
        \begin{theorem}\label{teo-ful}
Let   $t>0$, $\tilde x\in\R^n$ and
   \begin{equation}\label{sol form-0}
   u_{\tilde x, t}(x)=\beta_n \Big(\frac{t}{t^2+|x-\tilde x|^2}\Big)^{\frac n2},\quad\forall\, x\in \R^n,
    \end{equation}
where
 \begin{equation}\label{stand sol con}
 \beta_n=2^{\frac n2} e^{\frac n2\psi(\frac n2) }.
 \end{equation}
 Then problem \eqref{eq 1.1} only admits the solutions' set $\{u_{\tilde x, t}\}_{t>0,\,\tilde x\in\R^n}$.
    \end{theorem}

Here we mention that (\ref{eq 1.1}) and the solutions with the form (\ref{sol form-0}) are approached by the ones of critical fractional problem
 \begin{equation}\label{eq 1.1-s}
 (-\Delta)^\s  u(x)= u^{2^*_s-1}, \;\; u> 0 \;\;   {\rm in}\ \ \R^n,
\end{equation}
where  $n>2s$ and $2^*_s=\frac{2n}{n-2s}$.   It was proved in \cite[Theorem 1.1]{CLO}  (also see \cite{FL}) that all the solution of problem \eqref{eq 1.1-s} has the form of
 \begin{equation} \label{sol-s1}
 u_{s,t}(x)=b_{n,s} \Big(\frac{t}{t^2+|x-\tilde x|^2}\Big)^{\frac{n-2s}{2}},\ \ \forall\, x\in\R^n
 \end{equation}
for any $t>0$, where $ b_{n,s}>0$ is a constant (see (\ref{sharp const}) below) and $\tilde x\in\R^n$.  From the Caffarelli-Silvestre's extension \cite{Caffarelli-Silvestre},  equation (\ref{eq 1.1-s}) could be built by connection   with
the extension problem
\begin{equation}\label{eq 1.1-extension}
\left\{ \arraycolsep=1pt
\begin{array}{lll}
 \ \, {\rm div} (t^{1-2s} \nabla w) =  0    \quad \  &{\rm in}\quad  \R^{n+1}_+:=(0,+\infty)\times \R^n,\\[2mm]
 \phantom{    }
\displaystyle \lim_{t \to 0^+} t^{1-2s}\partial_t w  = w(0,\cdot)^{2^*_s-1} \quad \ &{\rm in}  \quad \R^n= \partial \R^{n+1}_+,
\end{array}
\right.
\end{equation}
via the boundary trace  $u=w_s(0,\cdot)$ in $\R^n$.  This type of extension  builds the platform for the Yamabe problem in the fractional setting, see the references  \cite{CG,CC,KMW,DSV}  and more related study could see  \cite{CK,JQS,NPX,CV2}.  An amazing extension involving the logarithmic Laplacian
has been built   \cite{CHW} recently.

 The explicit formula of $b_{n,s}$
is very important in our study of the approximating the solution of (\ref{eq 1.1}) by passing to the limit as $s\to0^+$. To this end,
we provide the proof of the sharp constant according to  \cite[Theorem 3.1]{Lieb}
 \begin{equation} \label{sharp const}
 b_{n,s}=2^{\frac{n-2s}{2}} \Big( \frac{\Gamma(\frac{n}{2}+s)}{\Gamma(\frac{n}{2}-s)}  \Big)^{\frac{n-2s}{4s}}.
  \end{equation}

   Our approach of the extremal solutions is to pass to limit of the ones of fractional Yamabe problems and
 to show the uniqueness by the method of Chen-Li-Ou in \cite{CLO} under the restriction of  $L^1_0(\R^n)\cap \cD^0(\R^n)$.
In fact, from our proof, we can also provide a correct sharp constant  of the Pitt's inequality by  Beckner \cite{Beckner}:   for $\|f\|_{L^2(\R^n)}=1$,
    \begin{align}\label{Pitt ineq}
    \frac n2 \int_{\R^n}(\ln|\xi|)|\hat{f}(\xi)|^2 d\xi \geq \int_{\R^n} (\ln |f|) |f|^2dx+B_n,
   \end{align}
 where
 \begin{align}\label{Pitt ineq-con}
 B_n=\ln{\|u_1\|_{L^2(\R^n)}}  =\frac n2\Big(\psi(\frac n2) +\frac12 \ln \pi -\frac1n \ln\big( \frac{\Gamma(n)}{\Gamma(\frac n2)} \big)  \Big).
    \end{align}
 We would like to mention that Beckner in \cite{Beckner1} (also see \cite{Beckner2}) shows that the Pitt's inequality in the unit sphere $\mathbb{S}^{n}$ of $\R^{n+1}$,
   \begin{align}\label{Pitt ineq-sph}
    \int_{\mathbb{S}^{n}\times \mathbb{S}^{n}} \frac{ \big|f(\xi)-f(\eta)\big|^2}{|\xi-\eta|^{n} } d\xi d\eta\geq  \frac2{\pi^{1/2}}\frac{\Gamma(\frac{n+1}{2})}{\Gamma(\frac{n}2 +1)}  \int_{\mathbb{S}^{n}}   \big( \ln |f|\big)|f|^2 dx,\quad  \forall\, f\in L^2(\mathbb{S}^n), \ \|f\|_{L^2(\mathbb{S}^n)}=1
   \end{align}
and  the extremal functions are of the form
 \begin{align}\label{Pitt ineq-sph-sol}
 A(1-\xi\cdot \eta)^{-\frac n2} \quad \text{ for    $\xi\in \mathbb{S}^n$  }
    \end{align}
 with $|\eta|<1$ and
 for some suitable constant $A$. Recently, Frank-K\"onig-Tang \cite{FKT}  gave the classification  of solutions for
 \begin{align}\label{Pitt ineq-con}
 \int_{\mathbb{S}^n} \frac{ u(\xi)-u(\eta)}{|\xi-\eta|^{n} } d\eta =E_n u(\xi) \ln u(\xi)   \quad {\rm for}\ \  \xi\in \mathbb{S}^n,
 \end{align}
  highlighted   the conformal invariance for the above equation and  showed that it admits the solutions with the form  of (\ref{Pitt ineq-sph-sol}),
  where $E_n=\frac4n \frac{\pi^{\frac n2}}{\Gamma(\frac n2)}$.
Notice that from the stereographic projection,  there is some 'equivalence' between our problem (\ref{eq 1.1}) and  problem (\ref{Pitt ineq-con})
under the suitable restriction, such as $u\in L^2(\R^n)$ for (\ref{Pitt ineq-con}).  Our approach enables us to classify the solution within the set of $L^1_0(\R^n)\cap \cD^0(\R^n)$.

 Here we remark that   the crucial point for classification lies in the method of moving planes.  Our challenges of moving planes arises  from the logarithmic Laplacian and non-monotonicity of nonlinearity. We cannot adopt the approach of  Chen-Li-Ou in \cite{CLO} via the integral equation,  as the fundamental
  solution of the logarithmic Laplacian   is no longer positive and decreasing-monotone. The core of the traditional method of moving planes is to involve the  Sobolev embedding to start the moving planes   from $-\infty$ in one direction, see e.g. \cite{FW}.  While  there is no  pertinent  Sobolev embeddings
  related to the logarithmic operator in $\R^n$ and   the nonlinearity is no longer monotonically increasing,
so it fails to apply  this type moving planes to give the classification.  Here we adopt a method of moving planes suitable for the attributes of the logarithmic Laplacian. 
Furthermore, another important tool for proving our uniqueness is  the Pitt's inequality (\ref{Pitt ineq})   and Beckner in \cite{Beckner} indicates that the equality in (\ref{Pitt ineq}) is achieved by  the   extremal functions  with the form $A(1+|x|^2)^{-\frac n2}$ via "the conformal invariance of this inequality is inherited from the Hardy-Littlewood-Sobolev inequality and this
 suffices to determine the extremals."

 The most significant features of (\ref{eq 1.1}) lie in the capability of employing the Fourier transform,  the scaling and Kelvin transform.
 These properties inherit from the logarithmic Laplacian and the logarithmic type nonlinearity and they  play an important role in procedure of the moving planes.

  \begin{proposition}\label{cr 3.0}
 Let  $T\in\R$,  $u\in \cD^0(\R^n)\cap L^1_0(\R^n)$ be a solution of
 \begin{equation}\label{eq 1.1-T}
\left\{ \arraycolsep=1pt
\begin{array}{lll}
\cL_\Delta u+Tu= \frac4n u\ln u \quad \  &{\rm in}\ \,   \R^n,\\[2mm]
 \phantom{ \cL_{\Delta} +Tu  }
  u\geq 0\quad \ &{\rm{in}}  \ \, \R^n.
\end{array}
\right.
\end{equation}

  $(i)$ Given $l>0$ and $x_0\in\R^n$, let
   \begin{align}\label{scal k-1}
  v_l(x)=l   u( x+x_0),\quad x\in\R^n.
   \end{align}
    Then $v_l$ verifies
    \begin{equation}\label{eq 1.1-l}
\left\{ \arraycolsep=1pt
\begin{array}{lll}
\cL_\Delta v_l + \big(T+\frac4n\ln l\big) v_l = \frac4n v_l\ln v_l\quad \  &{\rm in}\ \,  \R^n,\\[2mm]
 \phantom{ \cL_{\Delta} +T+ \big(\frac4n\ln l\big) uv_l   }
 v_l \geq 0\quad \ &{\rm{in}}  \ \,  \R^n.
\end{array}
\right.
\end{equation}

 $(ii)$     Given $l>0$, let
   \begin{align}\label{scal k-1}
  u_l(x)=l^{-\frac n2} u(l^{-1}x),\quad x\in\R^n.
   \end{align}
    Then $u_l$ is also a solution of (\ref{eq 1.1-T}).
Furthermore,  if $u\in L^2(\R^n)$, then
    $\|u_l\|_{L^2(\R^n)}= \|u \|_{L^2(\R^n)}$.
\smallskip

$(iii)$    Let $\tilde x\in\R^n, r>0$ and
    \begin{align}\label{def k-0}
    u^{\#_{0,r}}(x)  =\big( \frac{r}{|x- \tilde x|}\big)^{n} u(x^{*,r})\quad\ {\rm for}\ \, x\in\R^n\setminus\{\tilde x\},
      \end{align}
where
   \begin{align}\label{def k-1}
 x^{*,r}=\frac{r^2(x-\tilde x)}{|x-\tilde x|^2}+\tilde x \quad\ {\rm for}\ \, x\in\R^n\setminus\{\tilde x\}.
  \end{align}
 Then
  $$\cL_\Delta u^{\#_{0,r}}+ T u^{\#_{0,r}}=\frac4n u^{\#_{0,r}}\ln u^{\#_{0,r}} \quad\ {\rm in}\ \, \R^n\setminus\{\tilde x\}. $$
  \end{proposition}

    \begin{remark}


Let $\tilde u=e^{-\frac n4 \rho_n} u$, then $\tilde u$ verifies that
  \begin{equation}\label{eq 1.1-sd}
\left\{ \arraycolsep=1pt
\begin{array}{lll}
\cL_\Delta \tilde u-\rho_n \tilde u= \frac4n \tilde u\ln \tilde u \quad \  &{\rm in}\ \   \R^n,\\[2mm]
 \phantom{ \cL_{\Delta} ---  }
  \tilde u\geq 0\quad \ &{\rm{in}}  \ \ \R^n.
\end{array}
\right.
\end{equation}
This  form of equation (\ref{eq 1.1-sd})  is available for the method of  the moving planes.

    \end{remark}

Finally,  we show that the factor $\frac4n$ is the critical value for the existence of entire solution. In fact,  by a Pohozaev type identity,  we can derive  that    equation (\ref{eq 1.1}) has no positive solution, if $\frac4n$ is replaced by any other  positive constant.
  \begin{theorem}\label{teo 1.1-k}
Let $k>0$ and $u\in L^2(\R^n)\cap \cD^0(\R^n)$ be a nonnegative function verifying
$u^2 \ln u \in L^1(\R^n)$
and satisfying
\begin{equation}\label{eq 1.1-k}
\left\{ \arraycolsep=1pt
\begin{array}{lll}
\cL_\Delta u= k u\ln u \quad \  &{\rm in}\ \   \R^n,\\[2mm]
 \phantom{ \cL_{\Delta}   }
  u\geq 0\quad \ &{\rm{in}}  \ \  \R^n.
\end{array}
\right.
\end{equation}
If $k\not=\frac4n$, then $u\equiv 0$ in $\R^n$.
\end{theorem}

 Observe that for $k\not=\frac4n$, the Kelvin transform in Theorem \ref{cr 3.0} cannot be applied and we are unable to obtain the decay at  infinity
by the method of the moving planes. Hence, some restriction on $u$ need to be added to ensure the Pohozaev identity can derived and to get the desired result.

The rest of this paper is structured as follows. In Section 2, we demonstrate the qualitative properties of positive super solutions for  the critical logaritmic equation in Theorem \ref{cr 3.0}.    Section 3 is  dedicated to the existence of the positive solutions by approximating those of  fractional critical problems as the order $2s$ tends to zero. In Section 4  we adopt the method of moving planes to obtain the radially symmetric and asymptotic behavior at infinity and we perform the classification of positive solutions in Theorem \ref{teo-ful}. Finally, we show  the nonexistence of positive solutions of (\ref{eq 1.1-k}) in Theorem \ref{teo 1.1-k}

  \setcounter{equation}{0}
 \section{Preliminary }
 \subsection{Basic properties for Logarithmic Laplacian }

We start from the basic properties of  Eq. (\ref{eq 1.1}), including the scaling  and Kelvin transformation. \medskip

  \noindent{\bf Proof of Proposition \ref{cr 3.0}.}  $(i)$ For    a classical solution $u$ of (\ref{eq 1.1-T}),
  \begin{align*}
\cL_\Delta v_l   =  l \cL_\Delta u&= l \frac4n u\ln u - l Tu
  \\[1mm]&=\frac4n v_l \ln (v_l l^{-1})-Tv_l =\frac4n v_l \ln (v_l)-\big(T+\frac4 n\ln l\big)  v_l\quad {\rm in}\ \, \R^n.
  \end{align*}

\noindent $(ii)$  Recall  that  $u_l$ is defined in (\ref{scal k-1}).
   Note that
 $$\hat{ u}_l(\xi)=l^{\frac n2} \hat{u}(l\xi) $$
 and
  \begin{align*}
\cF(\cL_{\Delta}  u_l) (\xi)  =(2\ln|\xi|) \hat{ u_l}(\xi) &=   (2\ln|\xi|) l^{\frac n2}  \hat{ u} (l\xi)
  \\[1mm]&=l^{\frac n2}  (2\ln|l\xi|) \hat{ u} (l\xi)-l^{\frac n2}  (2\ln l) \hat{ u}  (l\xi),
  \end{align*}
  then it follows by the inverse Fourier transform
  \begin{align*}
 \cL_{\Delta}  u_l  (x)  &=l^{-\frac n2}  \big(\cL_{\Delta} u\big) (l^{-1}x)-l^{-\frac n2}  (2\ln l)u(l^{-1}x)
  \\[1mm]&= l^{-\frac n2}\Big( \frac4n u(l^{-1}x) \ln (u(l^{-1}x)) -Tu(l^{-1}x) \Big) -  (2\ln l)u_l(x)
  \\[1mm]&=    u_l(x) \big(\frac4n \ln (u_l(x))+\frac4n \ln l^{\frac{n}{2}}-T\big)  -  (2\ln l)u_l(x)
   \\[1mm]&=   \frac4n u_l(x)\ln (u_l(x)) - Tu_l(x),
  \end{align*}
which means that $u_l$ is  a solution of (\ref{eq 1.1-T}).

It is obvious that
 \begin{align*}
 \|  u_l\|^2_{L^2(\R^n)}    = {\int_{\R^n} |u(l^{-1} x)|^2 l^{-n} dx } =  {\int_{\R^n} |u(y)|^2  dy } ,
  \end{align*}
which is independent of $l>0$.  \smallskip

  $(iii)$  Given $r>0$ and a function $u\in C^\gamma(\R^n)\cap L^1(\R^n)$ with $\gamma>0$,
  let
   \begin{align}\label{def k-1}
 x^{*,r}=\big(\frac{r}{|x-\tilde x|^2}\big)^2(x-\tilde x)+\tilde x,  \quad u^{*,r}(x)=u(x^{*,r}),\quad  u^{\#_{s,r}}(x)=\big(\frac{r}{|x-\tilde x|}\big)^{n-2s} u^{*,r}(x)
  \end{align}
 for  $\tilde x\in\R^n$ and $x\in\R^n\setminus\{\tilde x\}$.
  Particularly, we use the notations that
    \begin{align}\label{def k-0}
   u^{\#_{0,r}}(x)=\big(\frac{r}{|x-\tilde x|}\big)^{n}  u^{*,r}(x)\quad\ {\rm for}\ \, x\in\R^n\setminus\{\tilde x\}.
      \end{align}

It is known that the Kelvin transform for fractional Laplacian
   \begin{align}\label{Kel-s}
  (-\Delta)^su^{\#_{s,r}} (x)=  \big(\frac{r}{|x-\tilde x|}\big) ^{2s+n} \big((-\Delta)^s u^{*,r} \big)(x)\quad  \ {\rm for\ any}\ \, x\in\R^n\setminus\{\tilde x\}.
  \end{align}

 We claim that
  \begin{align}\label{Kel-0}
 \cL_\Delta u^{\#_{0,r}} (x)=   \big(\frac{r}{|x-\tilde x|}\big) ^{n}\big( \cL_\Delta u\big)(x^{*,r})-4\big(\ln  |x-\tilde x|  -\ln r\big) u^{\#_{0,r}}(x) \quad \ {\rm for\ any}\ \, x\in\R^n\setminus\{\tilde x\}.
  \end{align}

 In fact,  we observe that
   \begin{align*}
  u^{\#_{s,r}} (x) &= \big(\frac{r}{|x-\tilde x|}\big) ^{n} u^{*,r}(x)\Big( 1+ 2s\big(\ln |x-\tilde x|-\ln r \big) +U_s(x)\Big),
   \end{align*}
where
$$\limsup_{s\to0^+}\Big(\| s^{-2}U_{s}(x) \|_{L^\infty(\R^n)} +\|s^{-2} U_{s}(x) \|_{C^\gamma(\R^n)}\Big) <+\infty.  $$
  We apply Theorem \ref{teo 3.1} below in Section \S3 to obtain that
  \begin{align*}
 (-\Delta)^s u^{\#_{s,r}} (x) &=(I+s\cL_\Delta +O(s^2)\Big) \bigg( \big(\frac{r}{|x-\tilde x|}\big) ^{n} u^{*,r}(x)\Big( 1+ 2s\big(\ln |x-\tilde x| -\ln r \big) +U_s(x)\Big)\bigg)
 \\[1mm]&=  u^{\#_{0,r}} (x)+s\Big( 2u^{\#_{0,r}} (x)\ln\frac{ |x-\tilde x|}{r}+ \cL_\Delta u^{\#_{0,r}} (x)\Big)+O(s^2)
 \end{align*}
and
\begin{align*}
  &\quad   \big(\frac{r}{|x-\tilde x|}\big) ^{2s+n} \big((-\Delta)^s u\big)(x^{*,r})
    \\[1mm]&= \big(\frac{r}{|x-\tilde x|}\big) ^{n} \Big( 1- 2s\big(\ln  |x-\tilde x|  -\ln r\big)   +O(s^2)\Big)  \big(\cI+ s\cL_\Delta+O(s^2)\big)  u^{*,r} (x)
  \\[1mm]&= u^{\#_{0,r}} (x)+s\big(\frac{r}{|x-\tilde x|}\big) ^{n}  \Big(\cL_\Delta u^{*,r} (x) -2\big(\ln |x-\tilde x| -\ln r \big) u^{*,r} (x) \Big)+O(s^2).
 \end{align*}
 From the identity (\ref{Kel-s}), we obtain that
  $$  2u^{\#_{0,r}} (x)(\ln |x-\tilde x|-\ln r)+ \cL_\Delta u^{\#_{0,r}} (x) =\big(\frac{r}{|x-\tilde x|}\big) ^{n}  \Big(\cL_\Delta u^{*,r} (x) -2\big(\ln |x-\tilde x| -\ln r \big) u^{*,r} (x) \Big),$$
  which implies (\ref{Kel-0}).  \smallskip

 Now it follows from (\ref{Kel-0}) that
  \begin{align*}
 \cL_\Delta u^{\#_{0,r}} (x)&=  \big(\frac{r}{|x-\tilde x|}\big) ^{n} \big( \cL_\Delta u\big)(x^{*,r})-4\big(\ln |x-\tilde x|-\ln r\big) u^{\#_{0,r}}(x)
 \\[1mm] & =\big(\frac{r}{|x-\tilde x|}\big) ^{n} \Big(\frac4n u^{*,r}\ln u^{0,r} -Tu^{*,r} \Big)-4\big(\ln |x-\tilde x|-\ln r\big) u^{\#_{0,r}}(x)
  \\[1mm]&= \frac4n u^{\#_{0,r}}(x)\ln\big( u^{\#_{0,r}}(x)\big)  -Tu^{\#_{0,r}}(x).
 \end{align*}
We complete the proof.   \hfill$\Box$\ms

 \subsection{Sharp constant $b_{n,s}$}

 \noindent{\bf Proof of \eqref{sharp const}.} Recall that by \cite[Theorem 1.1]{CLO}, the solution of problem \eqref{eq 1.1-s} has the form \eqref{sol-s1}, that is
$$u_{s,t}(x)=b_{n,s} \Big(\frac{t}{t^2+|x|^2}\Big)^{\frac{n-2s}{2}},\ \ \forall\, x\in\R^n,$$
for any $t>0$ and some suitable $b_{n,s}>0$.  Next we give the precise formula of the constant $b_{n,s}$.  To this end, we take $t=1$.
Let
$$v_0(x)=b_{n,s}^{-1}u_{s,1}(x)=\Big(\frac1{1+|x|^2}\Big)^{\frac{n-2s}{2}},\ \ \, \forall\, x\in\R^n, $$
then
   \begin{equation}\label{eq 1.1-s-1}
    (-\Delta)^\s v_0(x)=b_{n,s}^{2^*_s-2} v_0^{2^*_s-1}(x),\ \ \, \forall\, x\in\R^n .
  \end{equation}
It is know that the fundamental solution of $(-\Delta)^s$ in $\R^n$ has the following form
$$
\Phi_s(x)= \cP_0(s)  |x|^{2s-n}, \quad \forall\,  x\in\R^n\setminus\{0\},
$$
   where
    \begin{equation}\label{fund s.1}
 \cP_0(s):=   \pi^{-\frac{n}{2}}\, 4^{-s} \frac{\Gamma(\frac{n-2s}{2})}{\Gamma(s)}.
  \end{equation}
From (\ref{eq 1.1-s-1}), we have that
 $$v_0(x)=b_{n,s}^{2^*_s-2} \big(\Phi_s\ast v_0^{2^*_s-1}\big)(x) ,\ \ \forall\, x\in\R^n,$$
 or equivalently
    \begin{equation}\label{fund s.0}
    v_1^{\frac1{2^*_s-1}}(x)=b_{n,s}^{2^*_s-2} \big(\Phi_s\ast v_1 \big)(x) ,\ \ \forall\, x\in\R^n,
    \end{equation}
 where
 $ 
\displaystyle  v_1(x):=  v_0^{2^*_s-1}(x)=\big(1+|x|^2\big)^{-\frac{n+2s}{2}},\ \ \forall\, x\in\R^n.
   $ 

     From \cite[Theorem 3.1]{Lieb} with $\lambda=n-2s$, $p=\frac{2n}{n+2s}$ and $q=\frac{2n}{n-2s}$
     the best constant is defined as following and has the explicit formula:
     $$\cP_1(s):=\sup_{w\in L^q(\R^n) } \frac{\big\||x|^{2s-n} \ast w\big\|_{q}}{\|w\|_p}
     =\pi^{\frac{n-2s}{2}}\frac{\Gamma(s)}{\Gamma(\frac{n}{2}+s)}  \Big(\frac{\Gamma(\frac n2)}{\Gamma(n)}\Big)^{-\frac{2s}{n}}$$
and it is achieved uniquely by $v_1$, up to multiplications by nonzero constant and translations.

 Thus, we get
 \begin{align}\label{fund s.3}
  \sup_{w\in L^q(\R^n) } \frac{\big\|    \Phi_s \ast w\big\|_{q}}{\|w\|_p} &= \cP_1(s)  \cP_0(s)
 =(4  \pi)^{-s}\frac{\Gamma(\frac{n}{2}-s)}{\Gamma(\frac{n}{2}+s)}  \Big(\frac{\Gamma(\frac n2)}{\Gamma(n)}\Big)^{-\frac{2s}{n}}.
  \end{align}
On the other hand, by (\ref{fund s.0}) we have
 \begin{align}
  \sup_{w\in L^q(\R^n) } \frac{\big\|  \Phi_s \ast w\big\|_{q}}{\|w\|_p}
  &= \frac{ \Big(\displaystyle\int_{\R^n}  |\Phi_s \ast v_1|^{\frac{2n}{n-2s}}dx\Big)^{\frac{n-2s}{2n}} }{
  \Big(\displaystyle\int_{\R^n} v_1^{ \frac{2n}{n+2s} }dx\Big)^{\frac{n}{2n}}}\nonumber
   \\[1mm]& =b_{n,s}^{2-2^*_s}   \Big(\int_{\R^n} v_1^{  \frac{2n}{n+2s} }dx \Big)^{-\frac{2s}{n}}\nonumber
  \\[1mm]& =b_{n,s}^{2-2^*_s}   \Big(\int_{\R^n}   \big(1+|x|^2\big)^{-n} dx\Big)  ^{-\frac{2s}{n}}\label{sharp const1}
  \\[1mm]&=b_{n,s}^{2-2^*_s}  \pi^{-s} \Big(   \frac{ \Gamma(\frac n2)   }{\Gamma(n)} \Big)  ^{-\frac{2s}{n}},  \nonumber
    \end{align}
   where we have used
   $$ \int_{\R^n} \frac1{(1+|z|^2)^{\frac n2+\tau}}dz=\frac{\omega_{n-1}}2  B(\frac n2,\tau)=\frac{\omega_{n-1}}2  \frac{\Gamma(\frac n2)\Gamma(\tau)}{\Gamma(\frac n2+\tau)} $$
by taking $\tau =\frac{n}{2}$ and $\omega_{n-1}=|\mathbb S^{n-1}|=\frac{2\pi^{\frac n2}}{\Gamma(\frac n2)}.$
Compare with (\ref{fund s.3}) and (\ref{sharp const1}), we have that
     \begin{align*}
 b_{n,s}
 =\Big(2^{-2s}  \frac{\Gamma(\frac{n-2s}{2})}{\Gamma(\frac{n}{2})}\Big)^{-\frac{n-2s}{4s}}=2^{\frac{n-2s}{2}} \Big(\frac{\Gamma(\frac{n}{2}+s)}{\Gamma(\frac{n}{2}-s)}\Big)^{\frac{n-2s}{4s}},
    \end{align*}
that is  (\ref{sharp const}). \hfill$\Box$\medskip

 Next we consider the limits of the constant $b_{n,s}$ as $s\to0^+$.

   \begin{lemma}\label{lm 2.1}
   Let $b_{n,s}$ be defined in (\ref{sharp const}),
  $$b_{n,0}= e^{\frac n2\big( \ln2 +\psi(\frac{n}{2})\big)}\quad
  {\rm and}\quad  b_{n,1}= -\Big( \frac n2    \psi'(\frac{n}{2})+ \psi(\frac{n}{2})+\ln  2\Big) b_{n,0},$$
  where $\psi(\tau):=  \frac{d}{d\tau} \ln \Gamma(\tau)=\frac{\Gamma'(\tau)}{\Gamma(\tau)}$ is the digamma function.

 Then
   \begin{align}\label{lim s0}
(a) \qquad   \lim_{s\to0^+} b_{n,s}=b_{n,0} \;\quad \text{and} \;\quad (b) \qquad \lim_{s\to0^+} \frac{b_{n,s}-b_{n,0} }s=  b_{n,1}
    \end{align}
and
 \begin{align}\label{lim s2}
   b_{n,s}=b_{n,0}  +s  b_{n,1}+ O(s^2) \quad {\rm as}\ \,  s\to0^+.
    \end{align}

   \end{lemma}
 \noindent{\bf Proof. } A direct computation shows that
  \begin{align*}
  \lim_{s\to0^+} \ln  b_{n,s}&= \frac12 \lim_{s\to0^+} (n-2s)\Big(\frac{\ln  \Gamma(\frac{n}{2}+s)-\ln  \Gamma(\frac{n}{2}-s) }{2s}\Big)+ \ln  2\lim_{s\to0^+} \frac{n-2s}{2}
  \\[1mm]&=\frac n4  \lim_{s\to0^+}\Big(\psi(\frac{n}{2}+s)+\psi(\frac{n}{2}-s)\Big)+\frac n2 \ln 2
   \\[1mm]&= \frac n2\Big( \ln 2 +\psi(\frac{n}{2})\Big),
       \end{align*}
 which implies part $(a)$ of (\ref{lim s0}).  Note that
  \begin{align*}
 \lim_{s\to0^+}  \frac{d}{ds} (\ln  b_{n,s})&= \lim_{s\to0^+} \frac{d}{ds} \Big( \frac{n-2s}{4s}\big(\ln  \Gamma(\frac{n}{2}+s)-\ln  \Gamma(\frac{n}{2}-s) \big)+ \frac{n-2s}{2} \ln  2\Big)
   \\[1mm]&= \frac n4 \lim_{s\to0^+}\bigg[  \frac{1}{s}   \Big(\psi(\frac{n}{2}+s)+\psi(\frac{n}{2}-s)\Big) - \frac{1}{s^2}\big(\ln  \Gamma(\frac{n}{2}+s)-\ln  \Gamma(\frac{n}{2}-s) \big)\bigg]
  \\[1mm]&\qquad-\frac12 \lim_{s\to0^+} \Big(\psi(\frac{n}{2}+s)+\psi(\frac{n}{2}-s)\Big) - \ln 2
   \\[1mm]&= -\frac n2    \psi'(\frac{n}{2}) -  \psi(\frac{n}{2})-\ln  2,
       \end{align*}
 here $\psi'$ is the derivative of  $\psi$ and it is continuous.
 As a consequence, we obtain that
    \begin{align*}
 b_{n,1}:= \lim_{s\to0^+} \frac{b_{n,s}- b_{n,0}}s
  &=  \lim_{s\to0^+} b_{n,s} \frac{d}{ds}  ( \ln  b_{n,s})=-\Big( \frac n2    \psi'(\frac{n}{2})+ \psi(\frac{n}{2})+\ln  2\Big) b_{n,0},
       \end{align*}
  which implies part $(b)$ of (\ref{lim s0}).

  By similarly proof,  set
   \begin{align*}
b_{n,2}:=  \lim_{s\to0^+} \frac{b_{n,s}-b_{n,0}-b_{n,1}s}{s^2} \in(0,+\infty),
       \end{align*}
then
$$b_{n,s}=b_{n,0}+b_{n,1}s+O(s^2)\quad {\rm as}\ \, s\to0^+. $$
We complete the proof.  \hfill$\Box$\medskip

 \setcounter{equation}{0}
 \section{Existence of bubble solutions}

\subsection{Expansion for the fractional bubble solutions}

We have proved the expansion of the sharp constant $b_{s,t}$ as $s \to 0$ in the last subsection. We need to prove that the limit
of $u_{s,t}$, the solution of equation  \eqref{eq 1.1-s}, is in fact a solution of \eqref{eq 1.1}. More precisely we prove and then use the following convergence result which should be useful for other situation.

 \begin{theorem}\label{teo 3.1}
 Let $s_0\in(0,\frac14)$, $s\in(0,\frac18)$, $q\in(1,+\infty)$ and $\{u_s\}_{s\in(0,s_0)}\subset C^2_{\loc}(\R^n)\cap L^1_0(\R^n)$
 be a sequence of functions such that
 \begin{align}\label{assump-1}
 \limsup_{s\to0^+}\Big(\sup_{x\in\R^n} \big(|u_s (x)|+ |\nabla u_s (x)| \big)+ \big\| (|u_s|+ |\nabla u_s| )\big\|_{L^q(\R^n)}\Big)<+\infty.
  \end{align}
 Suppose that there exist $u_0,u_1\in C^2_{\loc}(\R^n)\cap L^1_0(\R^n)$ verifying
 \begin{align}\label{assump-2}
 \sup_{x\in\R^n}\big(|u_0(x)|+|u_1(x)|+ |\nabla u_0(x) | + |\nabla  u_1(x)|\big)+ \big\||u_0|+|u_1|+ |\nabla u_0 | + |\nabla  u_1|\big\|_{L^q(\R^n)} <+\infty
  \end{align}
and
 $$w_s:=\frac{1}{s^2} \big(u_s- u_0-su_1\big)\quad{\rm in}\ \R^n, $$
 satisfying
  \begin{align}\label{assump-3}
   \limsup_{s\to0^+}\Big( \sup_{x\in\R^n} \big( |w_s(x)|+  |\nabla w_s(x)|\big)+\big\| |w_s|+  |\nabla w_s|\big\|_{L^q(\R^n)} \Big)   <+\infty .
  \end{align}
Then as $s\to0^+$,
 $$(-\Delta)^s u_s = u_0+s\big(\cL_\Delta u_0 + u_1 \big)+o(s)\quad\text{in \ $L^p(\R^n)$ with $p\in(q,+\infty)$} $$
 and
 $$(-\Delta)^s u_s(x)= u_0(x)+s\big(\cL_\Delta u_0(x) + u_1(x) \big)+o(s)\quad\text{ uniformly for any $x\in\R^n$ }.  $$

 \end{theorem}
 \noindent{\bf Proof. } Let $u\in C^\beta_{\loc}(\R^n)\cap L^1_0(\R^n)$ with $\beta\in(0,1)$,  and for $0< \s < \min \{\frac{\beta}{2},\frac{1}{2}\}$, the principal value in the definition of $[(-\Delta)^\s u](x)$ reduces to a standard Lebesgue integral.
 For $x \in \R^n$, we then have
 \begin{align}
 & \frac1s\bigg((-\Delta)^\s u_s (x) -\Big(u_0(x)+s \big(\cL_\Delta u_0 (x)+u_1(x) \big) \Big)\bigg)\nonumber
 \\[8mm]&=\frac1s\bigg( \Big(c_{n,\s}\int_{\R^n \setminus B_1}|z|^{-n-2\s}dz \Big)\, u_s(x)- u_0(x)-s\big(\rho_nu_0(x)+ u_1(x)\big)\bigg)  \nonumber
 \\[1mm]&\quad -\frac1s\bigg(c_{n,\s}\int_{\R^n \setminus B_1}\frac{u_s(x+z)}{|z|^{n+2\s}}\,dz-c_{n}s\int_{\R^n \setminus B_1}\frac{u_0(x+z)}{|z|^{n}}\,dz\bigg)\label{mmao-1}
 \\[1mm]&\quad +\frac1s\bigg( c_{n,\s}\int_{B_1}\frac{u_s(x)-u_s(x+z)}{|z|^{n+2\s}}\,dz-c_ns \int_{B_1}\frac{u_0(x)-u_0(x+z)}{|z|^{n}}\,dz\bigg),\nonumber
 \end{align}
 where   $B_1 := B_1(0)$.

{\it Part 1:}   Let
$$d_{n}(\s) := \frac{c_{n,\s}}{\s}= 2^{2\s} \pi^{-\frac{n}{2}}
\frac{\Gamma(\frac{n}{2}+\s)}{\Gamma(1-\s)}, $$
and then
$$
d_n(0) = \pi^{-\frac{n}{2}}\Gamma(\frac{n}{2})= c_n \quad \text{and}\quad d_n'(0)= \pi^{-\frac{n}{2}} \Bigl( (2 \ln  2 - \gamma)\Gamma(\frac{n}{2})+ \Gamma'(\frac{n}{2})\Bigr)= c_n \rho_n.
$$
We now can write
\begin{equation}
  \label{eq:def-d-N-s}
c_{n,\s} = c_n  \big(1 + \rho_n s+  \bar d(s)\big) \s \qquad \text{for}\ \, s\in(0,\frac18),
\end{equation}
where $$\limsup_{s\to0^+} |\bar d(s)|s^{-2}<+\infty.$$
Observe that, by (\ref{eq:def-d-N-s}),
  \begin{align*}
   c_{n,\s}\int_{\R^n \setminus B_1}|z|^{-n-2\s}dz  =  \frac{c_{n,\s} |\mathbb{S}^{n-1}|}{2\s}
 &=1 +\rho_ns+\bar d(s),
   \end{align*}
   and by the definition of $w_s$ and (\ref{assump-2}),
we have then
  \begin{align*}
A_1(\s,x)&:= \frac1{s } \bigg(\Big(c_{n,\s}\int_{\R^n \setminus B_1}|z|^{-n-2\s}dz \Big)\, u_s(x)-   u_0(x)-s\big( \rho_nu_0(x)+  u_1(x)\big)\bigg)
\\[1mm]&= \frac1{s } \Big(\big(1  +\rho_ns+\bar d(s)\big)\big( u_0+su_1(x)+s^2 w_s (x)\big)-   u_0(x)-s\big(\rho_n  u_0(x)+  u_1(x)\big)\Big)
\\[1mm]&= \frac1{s }   \bar d(s) u_s (x)+s\big(1  +\rho_ns+\bar d(s)\big) w_s(x) ,
  \end{align*}
thus, for $q\in(1,+\infty)$
  \begin{align}  \label{part 1}
  \quad \lim_{s\to0^+} \|A_1(\s,\cdot)\|_{L^q(\R^n)}
 & \leq C \lim_{s\to0^+} \bigg( \bar d(s)s^{-1} \|u_s\|_{L^q(\R^n)}+s \|w_s\|_{L^q(\R^n)}
\bigg)
 =0
  \end{align}
  and
  \begin{align}  \label{part 1}
  \lim_{s\to0^+} \sup_{x\in\R^n} |A_1(\s,x)|
 & \leq C \lim_{s\to0^+} s\bigg( \bar d(s)s^{-2}\sup_{x\in\R^n}   |u_s(x)|+ s \sup_{x\in\R^n}   |w_s(x)|
\bigg)
 =0
  \end{align}
  by the assumptions of (\ref{assump-1}) and (\ref{assump-2}). \smallskip

 {\it Part 2:}   We see that
  \begin{align*}
A_2(\s,x): &=\frac1s\bigg( c_{n,\s}\int_{\R^n \setminus B_1}\frac{u_s(x+z)}{|z|^{n+2\s}}\,dz-c_n s \int_{\R^n \setminus B_1}\frac{u_0(x+z)}{|z|^{n}}\,dz\bigg)
\\[1mm]&= d_n(s)
  \int_{\R^n \setminus B_1}\frac{u_0(x+z)+su_1(x+z)+s^2 w_s(x+z)}{|z|^{n+2\s}}\,dz
-c_n \int_{\R^n \setminus B_1}\frac{u_0(x+z)}{|z|^{n}}\,dz
\\[1mm]&= c_n    \int_{\R^n \setminus B_1}u_0(x+z)\frac1{|z|^{n}} \big(\frac{1}{|z|^{2\s}}-1\big)\,dz
 +c_n  \big( \rho_n s+  \bar d(s)\big)  \int_{\R^n \setminus B_1}\frac{u_0(x+z)}{|z|^{n+2s}}dz
\\[1mm]&\quad+d_n(s)      \int_{\R^n \setminus B_1}\frac{ su_1(x+z)+s^2w_s(x+z)}{|z|^{n+2\s}}\,dz.
   \end{align*}

 Note that for any $\epsilon>0$, there exists $R=R(\epsilon)>1$ such that
   \begin{align*}
c_n    \Big|\int_{\R^n \setminus B_R}u_0(x+z)\frac1{|z|^{n}} \big(\frac{1}{|z|^{2\s}}-1\big)\,dz\Big|&\leq c_n   \int_{\R^n \setminus B_R}|u_0(x+z)|\frac1{|z|^{n}}  \,dz
\\[1mm]&\leq c_n   \Big(\int_{\R^n \setminus B_R} |u_0(x+z)|^q dz\Big)^{\frac1q}  \Big(\int_{\R^n \setminus B_R}  |z|^{-nq'} dz\Big)^{\frac1{q'}}
\\[1mm]&\leq  C \|u_0\|_{L^q(\R^n)} R^{-n+\frac{n}{q'}}
\\[1mm]&\leq \epsilon.
  \end{align*}
Moreover, since
\begin{align}
  \int_{\R^n \setminus B_R}|u_0(x+z)|\frac1{|z|^{n}}  \,dz&=\int_{\R^n }|u_0(x-z)| (1_{\R^N\setminus B_R} |z|^{-n})  \,dz \label{since}
  \\[1mm]&=|u_0|\ast (1_{\R^N\setminus B_R} |\cdot|^{-n}), \nonumber
    \end{align}
we have then
   \begin{align}
\Big(\int_{\R^n \setminus B_R} \Big|c_n     \int_{\R^n \setminus B_R}u_0(x+z)\frac1{|z|^{n}} \big(\frac{1}{|z|^{2\s}}-1\big)\,dz\Big|^p dx\Big)^{\frac{1}p}
&\leq \||u_0|\ast ((1_{\R^N\setminus B_R} |\cdot|^{-n}))  \|_{L^p (\R^n)} \label{since2}
\\[1mm]&\leq    \|u_0\|_{L^q(\R^n)}    \|1_{\R^N\setminus B_R} |\cdot|^{-n}  \|_{L^\tau(\R^n)} \nonumber
\\[1mm]&=  C \|u_0\|_{L^q(\R^n)}  R^{ -n+\frac{n}{\tau}} \nonumber
\\[1mm]&\leq \epsilon, \nonumber
  \end{align}
where we used  the Young's inequality for $p \in (q, +\infty) $
  and
  $\tau>1$ verifies that
  $$\frac1\tau=1+\frac1{p} -\frac1q.  $$

  Now for the fixed $R>0$ above,  we have that
 \begin{align*}
c_n   \Big|\int_{B_R \setminus B_1}u_0(x+z)\frac1{|z|^{n}} \big(\frac{1}{|z|^{2\s}}-1\big)\,dz\Big|&\leq c_n   \s  \int_{B_R \setminus B_1}\frac{|u_0(x+z)|}{|z|^{n}} (2\ln |z|) dz
\\[1mm]&\leq  c_n   \s (2\ln R)  \int_{B_R \setminus B_1}\frac{|u_0(x+z)|}{|z|^{n}}  dz.
  \end{align*}
 Similar to the above computations, we have that
    \begin{align*}
\int_{B_R \setminus B_1}\frac{|u_0(x+z)|}{|z|^{n}}  dz &\leq     C   \|u_0\|_{L^q(\R^n)} \big(1-R^{-nq'+n}\big)^{\frac1q'} \leq  C   \|u_0\|_{L^q(\R^n)}
  \end{align*}
 and
   \begin{align*}
\Big\|\int_{B_R \setminus B_1}\frac{|u_0(x+z)|}{|z|^{n}}  dz\Big\|_{L^p(\R^n)} & \leq    \|u_0\|_{L^q(\R^n)}  \big\|1_{B_R\setminus B_1} |\cdot|^{-n}  \big\|_{L^\tau(\R^n)}
\\[1mm]&\leq   C \|u_0\|_{L^q(\R^n)}
  \end{align*}
Consequently,
$$\lim_{s\to0^+} \sup_{x\in\R^n} \Big|  c_n \int_{\R^n \setminus B_1}u_0(x+z)\frac1{|z|^{n}} \Big(\frac{1}{|z|^{2\s}}-1\Big)\,dz\Big| =0 $$
and
$$\lim_{s\to0^+}    \Big\| c_n\int_{\R^n \setminus B_1}u_0(\cdot+z)\frac1{|z|^{n}} \Big(\frac{1}{|z|^{2\s}}-1\Big)\,dz\Big\|_{L^p(\R^n)} =0.  $$

Now we deal with the term $\int_{\R^n \setminus B_1}\frac{ su_1(x+z)+s^2 w_s(x+z)}{|z|^{n+2\s}}\,dz$. Note that
 \begin{align*}
d_n(s)    \Big| \int_{\R^n \setminus B_1}\frac{ su_1(x+z)+s^2 w_s(x+z)}{|z|^{n+2\s}}\,dz\Big| &
\leq Cs     \int_{\R^n \setminus B_1}\frac{|u_1(x+z)|+ s|w_s(x+z)|}{|z|^{n}} dz
\\[1mm]&\leq Cs   \Big( \|u_1\|_{L^q(\R^n)}+\|s w_s\|_{L^q(\R^n)}\Big) \Big(  \int_{\R^n \setminus B_1}\Big(    \frac{1}{|z|^{nq'}} dz\Big)^{\frac1{q'}}
\\[1mm]&\to0\quad {\rm uniformly \, in }\ \, x \ \, {\rm as}\ \, s\to0^+
  \end{align*}
  and
   \begin{align*}
d_n(s)    \Big\| \int_{\R^n \setminus B_1}\frac{ su_1(x+z)+s^2 w_s(x+z)}{|z|^{n+2\s}}\,dz\Big\|_{L^p(\R^n)}
 &\leq Cs    \Big\|  \int_{\R^n \setminus B_1}\frac{|u_1(x+z)|+ s|w_s(x+z)|}{|z|^{n}} dz\Big\|_{L^p(\R^n)}
\\[1mm]&\leq Cs   \Big( \|u_1\|_{L^q(\R^n)}+s \|  w_s\|_{L^q(\R^n)}\Big)
\\[1mm]&\to0\quad {\rm as}\ \, s\to0^+ ,
  \end{align*}
where we have used in the second inequality the Young's inequality for $p \in (q, +\infty) $ as in (\ref{since2}).
  Similar proof could be done for $c_n  \big( \rho_n s+  \bar d(s)\big)  \displaystyle\int_{\R^n \setminus B_1}\frac{u_0(x+z)}{|z|^{n+2s}}dz$.

Thus, we obtain that for $p\in(q,+\infty)$
 \begin{align} \label{part 2}
 \lim_{s\to0^+} \sup_{x\in\R^n}|A_2(\s,x)|=0 \quad{\rm and} \quad  \lim_{s\to0^+} \|A_2(\s,\cdot)\|_{L^p(\R^n)}=0.
 \end{align}

 {\it Part 3:}   we deal with the last term
 \begin{align*}
A_3(\s,x): &=\frac1s\bigg( c_{n,\s}\int_{B_1}\frac{u_s(x)-u_s(x+z)}{|z|^{n+2\s}}\,dz-c_ns \int_{B_1}\frac{u_0(x)-u_0(x+z)}{|z|^{n}}\,dz\bigg)
\\[1mm]&= d_n(s)  \int_{B_1}\frac{u_0(x)-u_0(x+z)+s\big(u_1(x)-u_1(x+z)\big) }{|z|^{n+2\s}}\,dz
\\[1mm]&\quad+ d_n(s)  \int_{B_1}\frac{ s^2\big(w_s(x)-w_s(x+z)\big)}{|z|^{n+2\s}}\,dz
 -c_n \int_{ B_1}\frac{u_0(x)-u_0(x+z)}{|z|^{n}}\,dz
\\[1mm]&= c_n    \int_{  B_1}\big(u_0(x)-u_0(x+z)\big)\Big(\frac{1}{|z|^{n+2\s}}-\frac1{|z|^{n}} \Big)\,dz
\\[1mm]&\quad +c_n  \big( \rho_n s+  \bar d(s)\big)   \int_{  B_1}\frac{u_0(x)-u_0(x+z)}{|z|^{n+2s}} dz
\\[1mm]&\quad +d_n(s)      \int_{  B_1}\frac{ s\big(u_1(x)- u_1(x+z)\big)+ s^2 \big(w_s(x)- w_s(x+z)\big)}{|z|^{n+2\s}}\,dz.
   \end{align*}

 Note that for any $\epsilon>0$, there exists $\delta_\epsilon\in(0,1)$ such that for $s \in (0,s_0)$
   \begin{align*}
    \Big|\int_{B_{\delta_\epsilon}}\big(u_0(x)-u_0(x+z)\big)\Big(\frac{1}{|z|^{n+2\s}}-\frac1{|z|^{n}} \Big)\,dz\Big|&\leq
    2  \int_{  B_{\delta_\epsilon}}\big|u_0(x)-u_0(x+z)\big| \frac{1}{|z|^{n+2\s}} \,dz
\\[1mm]&\leq  C   |\nabla u_0(x)|  \int_{  B_{\delta_\epsilon}}  \frac{1}{|z|^{n+2\s-1}} \,dz
\\[1mm]&\leq   C |\nabla u_0(x)|\, \delta_{\epsilon}^{1-2s}
  \end{align*}
which implies that
for $q>1$
   \begin{align*}
 \Big( \int_{\R^n} \Big|\int_{  B_{\delta_\epsilon}}\big(u_0(x)-u_0(x+z)\big)\Big(\frac{1}{|z|^{n+2\s}}-\frac1{|z|^{n}} \Big)\,dz\Big|^qdx\Big)^{\frac1q}
 \leq    C  \|\nabla u_0\|_{L^q(\R^n)} \delta_{\epsilon}^{1-2s} \leq \epsilon.
  \end{align*}
  Similarly, we have
  \begin{align*}
  \Big|\int_{  B_{\delta_\epsilon}}\frac{ s\big(u_1(x)- u_1(x+z)\big)+s^2\big(w_s(x)- w_s(x+z)\big)}{|z|^{n+2\s}}dz\Big|
  \leq C \big( s |\nabla u_1(x)|+ s^2 |\nabla w_s(x)|\big) \delta_{\epsilon}^{1-2s},
  \end{align*}
and
    \begin{align*}
  \Big( \int_{\R^n}  \Big|\int_{  B_{\delta_\epsilon}} & \frac{ s\big(u_1(x)- u_1(x+z)\big)+s^2\big(w_s(x)- w_s(x+z)\big)}{|z|^{n+2\s}}dz\Big|^qdx\Big)^{\frac1q}
 \\[1mm]&\leq   C  \Big( s\|\nabla u_1\|_{L^q(\R^n)}+ s^2 \|\nabla w_s\|_{L^q(\R^n)}\Big) \delta_{\epsilon}^{1-2s}
\\[1mm]&\leq    \epsilon.
  \end{align*}

 For fixed $\delta_\epsilon$,  it follows
 that there exists $s_1 >0$ small enough such that $\forall \, s \in (0,s_1)$,
  \begin{align*}
    \Big|\int_{   B_1\setminus B_{\delta_\epsilon} }\big(u_0(x)-u_0(x+z)\big)& \Big(\frac{1}{|z|^{n+2\s}}-\frac1{|z|^{n}} \Big)\,dz\Big|
 \\[1mm]& \leq s C |\nabla u_0(x)|     \int_{   B_1\setminus B_{\delta_\epsilon} }   \frac{1}{|z|^{n-1  }} \ln \frac{1}{|z|} \,dz
 \\[1mm]& \leq s C  |\nabla u_0(x)|  \frac{1}{\delta_\epsilon^{n-1  }} \ln \frac{1}{\delta_\epsilon} \leq \epsilon,
  \end{align*}
  where we used the following fact: the inequality $\frac{1}{r^{2s}} -1 \leq \ln\frac{1}{r^{4s}}$ holds for small $s >0$ and $ \forall \; r \in (\delta_\epsilon,1)$.
  Moreover it follows that for $q\in(1,+\infty)$,
   \begin{align*}
    \Big\|\int_{   B_1\setminus B_{\delta_\epsilon} }\big(u_0(x)-u_0(x+z)\big)& \Big(\frac{1}{|z|^{n+2\s}}-\frac1{|z|^{n}} \Big)\,dz\Big\|_{L^q(\R^n)}
 \\[1mm]& \leq s C \|\nabla u_0(x)\|_{L^q(\R^n)}    \frac{1}{\delta_\epsilon^{n-1  }} \ln \frac{1}{\delta_\epsilon} \leq \epsilon.
  \end{align*}
Furthermore, we have that
   \begin{align*}
 \Big|c_n  \big( \rho_n s+  \bar d(s)\big)  \int_{  B_1}\frac{u_0(x)-u_0(x+z)}{|z|^{n+2s}} dz \Big|
 \leq Cs |\nabla u_0(x)|\int_{B_1} |z|^{1-2s-n}dz \leq Cs |\nabla u_0(x)|
      \end{align*}
      and
     \begin{align*}
 \Big|    \int_{   B_1  }   \frac{ s\big(u_1(x)- u_1(x+z)\big)+s^2\big(w_s(x)- w_s(x+z)\big)}{|z|^{n+2\s}}\,dz \Big|
   \leq C s \Big(  |\nabla u_1(x)|+ s |\nabla w_s(x)|\Big)
        \end{align*}
The left proof is similar.

 Consequently, we derive that
 \begin{align} \label{part 3}
 \lim_{s\to0^+} \sup_{x\in\R^n}|A_3(\s,x)|=0 \quad{\rm and} \quad  \lim_{s\to0^+} \|A_3(\s,\cdot)\|_{L^p(\R^n)}=0\ \ \text{for $q\in(1,+\infty)$}.
 \end{align}

Dragging  the estimates (\ref{part 1}), (\ref{part 2}) and   (\ref{part 3}) into (\ref{mmao-1}), we derive that
$$
\lim_{\s \to 0^+}\sup_{x\in\R^n}\Big| \frac{(-\Delta)^\s u_s(x) -u_0(x)}{\s}-u_1(x)- \cL_\Delta u_0(x)  \Big|=0$$
and
$$
\lim_{\s \to 0^+}
\Big\|\frac{(-\Delta)^\s u_s -u_0}{\s}-u_1- \cL_\Delta u_0 \Big\|_{L^p(\R^n)}=0\quad\ \text{for every $p \in (q,\infty)$.}
$$
We complete the proof.   \hfill$\Box$\smallskip

\subsection{Approximation of solutions}
 Next we show the existence of positive solution of \eqref{eq 1.1}.
          \begin{theorem}\label{teo0}
Let   $t>0$ and
  \begin{equation}\label{sol form-00}
   u_t(x)=\beta_n \Big(\frac{t}{t^2+|x|^2}\Big)^{\frac n2},\quad\forall\, x\in \R^n,
     \end{equation}
where we recall
 $$
 \beta_n=2^{\frac n2} e^{\frac n2\psi(\frac n2) }.
$$
    Then $u_t$  is a  solution  of \eqref{eq 1.1}.
    \end{theorem}

 \noindent{\bf Proof.} Recall that $\{u_{s,t}\}$ given by \eqref{sol-s1} are solutions of the equation \eqref{eq 1.1-s} with $$ b_{n,s}=b_{n,0}+s b_{n,1}+O(s^2), $$
 by Lemma \ref{lm 2.1}.

Letting
$$v_{s,t}(x)=\Big(\frac{t}{t^2+|x|^2}\Big)^{\frac{n-2s}{2}},\qquad v_{t}(x)=\Big(\frac{t}{t^2+|x|^2}\Big)^{\frac{n}{2}} \quad {\rm for}\ \, x\in\R^n,  $$
and
$$v_{s,t}=v_t-\frac2n sv_t\ln v_t+s^2 V_{s,t} \quad  {\rm in}\ \,    \R^n,$$
where by the Taylor's expansion, we have that
$$V_{s,t}(x)=v_{s',t}(x)(\ln v_{s',t}(x))^2\quad {\rm for\ some\ } s'\in[0,s].$$
This implies that
$$\limsup_{s\to0^+}\big(\sup_{x\in\R^n}  |V_{s,t}(x)| +\||V_{s,t}   \|_{L^q(\R^n)}\big) <+\infty\quad {\rm for\ any} \ q>1.  $$
Moreover, we have that
$$\partial_i V_{s,t}=\frac1{s^2} \big( \partial_i v_{s,t} +(\frac2n s-1) \partial_i v_t+ \frac2n s (\partial_i v_t) \ln v_t \big) \quad  {\rm in}\ \,    \R^n,$$
then
$$\limsup_{s\to0^+}\Big(\sup_{x\in\R^n} \big|\nabla V_{s,t}(x)\big|   +\|| \nabla V_{s,t} \|_{L^q(\R^n)}\Big) <+\infty\quad {\rm for\ any} \ q>1.  $$
Therefore
 \begin{align}\label{limteq-0}
 u_{s,t} =b_{n,0} v_{t} +s\big(b_{n,1}v_t-\frac{2}{n}b_{n,0} v_t\ln v_t    \big)+s^2  U_{s,t},
 \end{align}
where
$$\limsup_{s\to0^+}\Big(\sup_{x\in\R^n} \big(|U_{s,t}(x)|+|\nabla U_{s,t}(x)|\big)   +\||\big(|U_{s,t}|+|\nabla U_{s,t}|\big)   \|_{L^q(\R^n)}\Big) <+\infty\quad {\rm for\ any} \ q>1.  $$From (\ref{limteq-0}), we have that
  \begin{align}
  u_{s,t}^{\frac{n+2s}{n-2s}} &=\Big(b_{n,0} v_{t} +s\big(b_{n,1}v_t-\frac{2}{n}b_{n,0} v_t\ln v_t    \big)+  O(s^2) \Big)^{1+\frac{4s}{n-2s}} \nonumber
  \\[1mm]&=\big(b_{n,0} v_{t} +s\big(b_{n,1}v_t-\frac{2}{n}b_{n,0} v_t\ln v_t    \big)+  O(s^2) \big)\nonumber
  \\[1mm]&\quad \cdot \Big(1+\frac{4s}{n-2s}\ln \big(b_{n,0} v_{t} +s\big(b_{n,1}v_t-\frac{2}{n}b_{n,0} v_t\ln v_t    \big)+  O(s^2) \big)\Big)  \nonumber
   \\[1mm]&=  b_{n,0} v_{t} +s\big(b_{n,1}v_t-\frac{2}{n}b_{n,0} v_t\ln v_t    \big)+\frac{4s}{n-2s}  b_{n,0} v_{t} \ln (b_{n,0} v_{t}) +  O(s^2)\nonumber
  \\[1mm]&=  b_{n,0} v_{t} +s\big(b_{n,1}v_t-\frac{2}{n}b_{n,0} v_t\ln v_t  +\frac{4}{n}  b_{n,0} v_{t} \ln (b_{n,0} v_{t})  \big)+O(s^2).\label{limteq-1}
 \end{align}

We apply Theorem \ref{teo 3.1} with $u_s=u_{s,t}$, $u_0=b_{n,0} v_{t}$ and $u_1=b_{n,1}v_t(x)-\frac{2}{n}b_{n,0} v_t \ln v_t $ to obtain that
  \begin{align}
 (-\Delta)^s u_{s,t}(x) &=(\cI+s\cL_\Delta +O(s^2)\Big) \bigg( b_{n,0} v_{t}(x) +s\Big(b_{n,1}v_t(x)-\frac{2}{n}b_{n,0} v_t(x)\ln\big( v_t(x) \big)   \Big)+  U_{s,t}(x)\bigg) \nonumber
 \\[1mm]&= b_{n,0} v_{t}(x) +s\Big(b_{n,0} \cL_\Delta v_{t}(x) -\frac{2}{n}b_{n,0} v_t(x)\ln\big( v_t(x) +b_{n,1}v_t(x) \Big)+O(s^2).\label{limteq-2}
 \end{align}
 Note that
 $$(-\Delta)^s u_{s,t}(x)=u_{s,t}^{\frac{n+2s}{n-2s}}\quad{\rm in}\ \, \R^n, $$
 which, combining (\ref{limteq-1}) and (\ref{limteq-2}),  implies that
   \begin{align*}
   b_{n,0} \cL_\Delta v_{t} -\frac{2}{n}b_{n,0}\,  v_t \ln  v_t  +b_{n,1}v_t =b_{n,1}v_t-\frac{2}{n}b_{n,0}\, v_t\ln v_t  +\frac4n b_{n,0}\,  v_{t} \ln (b_{n,0} v_{t}),
  \end{align*}
which implies that
$$\cL_\Delta (b_{n,0} v_{t}) = \frac4n b_{n,0}\,  v_{t} \ln (b_{n,0} v_{t}).$$
Thus $ u_t= \beta_n v_{t}$ is the solution of (\ref{eq 1.1}) with $\beta_n :=b_{n,0}$.  We complete the proof. \hfill$\Box$


  \setcounter{equation}{0}

  \section{Radial symmetry}

  In this section, we consider the radial symmetry properties of the positive solution of
  \begin{equation}\label{eq 1.1sym}
\left\{ \arraycolsep=1pt
\begin{array}{lll}
\cL_\Delta v -\rho_n v= \frac4n v\ln v  \quad \  &{\rm in}\ \   \R^n,\\[2mm]
 \phantom{ \cL_{\Delta}   -\rho_n u }
  v\geq 0\quad \ &{\rm{in}}  \ \ \R^n.
\end{array}
\right.
\end{equation}

\begin{theorem}\label{teo 1sym}
Let  $n\geq 2$,  $v\in L^1_0(\R^n)\cap \cD^0(\R^n)$ be  a positive  solution of (\ref{eq 1.1sym}), then  $v$ is radially symmetric and monotonically decreasing about some point.

Furthermore, there exists $v_\infty\in(0,+\infty)$ such that
 \begin{equation}\label{con h}
 \lim_{|x|\to+\infty}|x|^{n} v(x)=v_\infty.
 \end{equation}
\end{theorem}

We employ the method of moving planes to establish the radially symmetric property.

 \subsection{Tools for the moving plane method}
To proceed with the moving plane method, we need the following notations: for $\lambda\in\R$,
$$
T_\lambda :=\left\{x \in \mathbb{R}^{n} \mid x_{1}=\lambda\right\},
$$
the reflection of the point $x \in \R^n$ about $T_\lambda$  is
$$
 x^\lambda :=\left(2 \lambda-x_{1}, x_{2}, \ldots, x_{n}\right),
$$
 and
$$
 H_\lambda:=\big\{x \in \mathbb{R}^{n} \mid x_{1}<\lambda\big\},\quad  \quad \tilde{H}_\lambda :=\big\{x \mid  x^\lambda \in  H_\lambda\big\}.
$$

 \begin{lemma}   \label{strong}
 Let $\lambda\in\R$, $\Omega$ be a  domain in $ H_\lambda$   and $w\in
\cD^0\left( \overline{\Omega} \right) \cap L_{0}^1 \left( \R^n \right)$ verify that
 \begin{equation}\label{eq 1.1-W}
 \left\{ \arraycolsep=1pt
\begin{array}{lll}
\mathcal{L}_{\Delta}w+W(x)w(x)\geq{0} \quad \
&{ \rm in}\ \Omega,
\\[2mm]
w(x)\geq{0}\quad \
&{ \rm in}\  H_\lambda \setminus \Omega,
\\[2mm]
w(x^{\lambda})=-w(x) \quad \
&{ \rm in}\  H_\lambda ,
\end{array}
\right.
\end{equation}
where  $W:\Omega\rightarrow\R$ is a given the function.

 If $w \geq 0$ in $\Omega$, then either $w\equiv0$ in $H_\lambda$ or $w > 0$ in $\Omega$.
 \end{lemma}
\noindent{\bf Proof. } Suppose that there is some point $x_0 \in H_\lambda$ such that $w (x_0)=0$, then
\begin{align*}
0  \leq \frac1{ c_{n}}\Big(\cL_\Delta  w(x_0)   +W(x_0)  w(x_0)  \Big)  &=    \int_{H_\lambda} \frac{-w(y)}{\left|x_{0}-y\right|^{n}}   d y+\int_{\tilde{H}_\lambda} \frac{-w (y)}{\left|x_{0}-y\right|^{n}}   d y
\\[2mm]&=    \int_{H_\lambda} \frac{-w (y)}{\left|x_{0}-y\right|^{n}}   d y+\int_{H_\lambda} \frac{-w (y^\lambda)}{\left|x_{0}-y^\lambda\right|^{n}}   d y
\\[2mm]&=       \int_{H_\lambda} \frac{-w (y) }{\left|x_{0}-y \right|^{n}}  d y+\int_{H_\lambda} \frac{w ({y})}{\left|x_{0}-y^\lambda\right|^{n}}  d y
\\[2mm]&=       \int_{H_\lambda}\Big( \frac{1}{\left|x_{0}-y^\lambda\right|^{n}} -\frac{1 }{\left|x_{0}-y \right|^{n}}   \Big)  w (y) d y
 \\[2mm]&\leq  0,
\end{align*}
 where the equality holds only for
$$
w \equiv 0\ \,  {\rm  \ in }\ \, H_\lambda
$$
since there holds $ \frac{1}{\left|x_{0}-y^\lambda\right|^{n}} -\frac{1 }{\left|x_{0}-y \right|^{n}} <0$ for $x_0, \ y\in H_\lambda$.
 \hfill$\Box$\medskip

Now we provide the maximum principle in the narrow region. Let
$$\bP_{a,b}=\big\{(x_1,x')\in(a-b,a+b)\times \R^{n-1}\big\}\ \  {\rm for}\ a\in\R,\ b>0.  $$

 \begin{lemma}   \label{Narrow A}   Let $\Omega$ be a  domain in $ H_\lambda$, $\lambda\in\R$,  $\epsilon>0$ such that
 $$ \Omega  \cap  \bP_{\lambda,\epsilon} \not=\emptyset $$
 and $w\in
\cD^0\left( \overline{\Omega} \right) \cap L_{0}^1 \left( \R^n \right)$ be a  function  such  that $w\not\equiv0$ in $H_\lambda$
and
 \begin{equation}\label{eq 3.1-R}
 \left\{ \arraycolsep=1pt
\begin{array}{lll}
\mathcal{L}_{\Delta}w+V(x)w(x)\geq{0} \quad \
&{ \rm in}\ \Omega,
 \\[2mm]
  \ \ w(x)\geq{0}\quad \ &{ \rm in}\   H_\lambda \setminus  \Omega\\[2mm]
w(x^{\lambda})=-w(x) \quad \
&{ \rm in}\  H_\lambda ,
\end{array}
\right.
\end{equation}
where the function $V:\Omega\rightarrow\R$ is lower bounded.
If $\Omega  \cap  \bP_{\lambda,\epsilon}$ is unbounded, we assume more that
 \begin{align}\label{decay inf}
 \liminf_{x\in \Omega  \cap  \bP_{\lambda,\epsilon},|x|\to+\infty} w(x)\geq0.
 \end{align}
Then for $\epsilon>0$ sufficiently small,  there is no $x_0\in \Omega\cap \bP_{\lambda,\epsilon} $  satisfying
$$ w(x_0)=  \inf_{x\in\Omega} w(x)<0.$$
That is, if the minimal of $w$ in $\Omega$ is negative, then it can't be achieved
  in  $\Omega\cap    \bP_{\lambda,\epsilon} $.

\end{lemma}
\noindent{\bf Proof. } By contradiction, we assume that there is  $x_{\lambda,\epsilon}\in\Omega\cap    \bP_{\lambda,\epsilon} $ such that
$$w(x_{\lambda,\epsilon})=\inf_{x\in \Omega} w(x)<0. $$
When $\Omega \cap  \bP_{\lambda,\epsilon} $ is unbounded,  (\ref{decay inf}) ensures that the  point achieving the  negative minimal is bounded.

Observe that
 \begin{align}\label{aa-00-1}
 \frac1{c_n} \cL_\Delta  w(x)   =  \cL_1w(x)+\cL_2w(x), \quad\forall\,  x\in\Omega_\lambda,
 \end{align}
where
$$\cL_1w(x)=  \int_{\R^n} \frac{w (x )-w (y)}{ |x -y|^{n}} 1_{B_1(x)}(y)  d y,\qquad \cL_2w(x)=\frac{\rho_n}{c_n} w(x)  - \int_{\R^n}\frac{w (y)}{|x -y|^{n}}1_{B_1^c(x)}(y) dy.
$$

Note that
  \begin{align}\label{lim-oo-1}
  \int_{H_{\lambda}} \frac{ 1_{B_1(x_{\lambda,\epsilon})}(y )+1_{B_1(x_{\lambda,\epsilon})}(y^{\lambda}) }{\left|x_{\lambda,\epsilon}-y^{\lambda}\right|^{n}  } dy\to +\infty\quad {\rm as}\ \epsilon\to+\infty.
  \end{align}
Thus, direct computations shows that
  \begin{align}
 \nonumber   \cL_{1}w(x_{\lambda,\epsilon})  &=
 \int_{H_{\lambda}} \frac{w \left(x_{\lambda,\epsilon}\right)-w (y)}{\left|x_{\lambda,\epsilon}-y\right|^{n}}  1_{B_1(x_{\lambda,\epsilon})}(y) d y+\int_{\tilde{H}_{\lambda}} \frac{w \left(x_{\lambda,\epsilon}\right) -w (y)}{\left|x_{\lambda,\epsilon}-y\right|^{n}} 1_{B_1(x_{\lambda,\epsilon})}(y) d y
\nonumber \\[1mm]&=  \int_{H_{\lambda}} \frac{w \left(x_{\lambda,\epsilon}\right) -w (y)}{\left|x_{\lambda,\epsilon}-y\right|^{n}} 1_{B_1(x_{\lambda,\epsilon})}(y) d y+\int_{ H_\lambda} \frac{w \left(x_{\lambda,\epsilon}\right)-w (y^{\lambda})}{\left|x_{\lambda,\epsilon}-y^{\lambda}\right|^{n}}   1_{B_1(x_{\lambda,\epsilon})}(y^{\lambda}) d y
\nonumber  \\[1mm]&<      \int_{H_{\lambda}} \frac{w \left(x_{\lambda,\epsilon}\right)-w (y)}{\left|x_{\lambda,\epsilon}-y^{\lambda}\right|^{n}}  1_{B_1(x_{\lambda,\epsilon})}(y) d y+\int_{H_{\lambda}} \frac{w \left(x_{\lambda,\epsilon}\right)+w (y)}{\left|x_{\lambda,\epsilon}-y^{\lambda}\right|^{n}} 1_{B_1(x_{\lambda,\epsilon})}(y^{\lambda})  d y
\nonumber \\[1mm]&=   w  (x_{\lambda,\epsilon})  \int_{H_{\lambda}} \frac{ 1_{B_1(x_{\lambda,\epsilon})}(y )+1_{B_1(x_{\lambda,\epsilon})}(y^{\lambda}) }{\left|x_{\lambda,\epsilon}-y^{\lambda}\right|^{n}  } dy\label{pp-1-1}
\end{align}
 and
 \begin{align*}
 \cL_{2}w(x_{\lambda,\epsilon})&= \frac{\rho_n}{c_n} w (x_{\lambda,\epsilon}) - \int_{H_{\lambda}} \frac{ w (y)}{\left|x_{\lambda,\epsilon}-y\right|^{n}}  1_{B_1^c(x_{\lambda,\epsilon})}(y) d y
 -\int_{\tilde{H}_{\lambda}} \frac{ w (y)}{\left|x_{\lambda,\epsilon}-y\right|^{n}} 1_{B_1^c(x_{\lambda,\epsilon})}(y) d y
 \\[1mm]&=\frac{\rho_n}{c_n} w (x_{\lambda,\epsilon})   - \int_{H_{\lambda}} \frac{ w (y)}{\left|x_{\lambda,\epsilon}-y\right|^{n}}  1_{B_1^c(x_{\lambda,\epsilon})}(y) d y+\int_{ H_{\lambda}} \frac{ w (y)}{\left|x_{\lambda,\epsilon}-y^{\lambda} \right|^{n}} 1_{B_1^c(x_{0,\lambda})}(y^\lambda) d y
 \\[1mm]&=\frac{\rho_n}{c_n} w (x_{\lambda,\epsilon})  - \int_{H_{\lambda}} w (y) 1_{B_1^c(x_{\lambda,\epsilon})}(y)   \Big( \frac{1}{\left|x_{\lambda,\epsilon}-y\right|^{n}}- \frac{1}{\left|x_{\lambda,\epsilon}-y^{\lambda} \right|^{n}}\Big) dy
 \\[1mm]&\qquad \qquad\  -\int_{H_{\lambda}} w (y)  \frac{1_{B_1^c(x_{\lambda,\epsilon})}(y)-1_{B_1^c(x_{\lambda,\epsilon})}(y^{\lambda})   }{\left|x_{\lambda,\epsilon}-y^{\lambda}\right|^{n}}   dy,
\end{align*}
where $x_{\lambda,\epsilon}=((x_{\lambda,\epsilon})_1,x_{\lambda,\epsilon}')\in \Omega \cap \bP_{\lambda,\epsilon}$ and
 \begin{align*}
 - \int_{H_{\lambda}} & w (y) 1_{B_1^c(x_{\lambda,\epsilon})}(y)   \Big( \frac{1}{\left|x_{\lambda,\epsilon}-y\right|^{n}}- \frac{1}{\left|x_{\lambda,\epsilon}-y^{\lambda} \right|^{n}}\Big) dy
\\[1mm] &\leq-  w (x_{\lambda,\epsilon})\Big( \int_{B_1(0)^c}\frac{1 }{ |z|^{n+1}}      dz\Big)\,  |x_{\lambda,\epsilon} -x_{\lambda,\epsilon}^{\lambda,\epsilon}|
\\[1mm] &\leq -C|(x_{\lambda,\epsilon})_1-\lambda_{\lambda,\epsilon}| \,  w (x_{\lambda,\epsilon})
\\[1mm] &\leq -C \epsilon  \,  w (x_{\lambda,\epsilon})
\end{align*}
 and
  \begin{align*}
 -\int_{H_{\lambda}} w (y)  \frac{1_{B_1^c(x_{\lambda,\epsilon})}(y)-1_{B_1^c(x_{\lambda,\epsilon})}(y^{\lambda})   }{\left|x_{\lambda,\epsilon}-y^{\lambda}\right|^{n}}   dy
 \leq -C\epsilon\,  w (x_{\lambda,\epsilon}).
\end{align*}
As a consequence, we obtain that
\begin{align*}
 0&\leq \cL_\Delta   w (x_{\lambda,\epsilon}) + V(x_{\lambda,\epsilon})   w (x_{\lambda,\epsilon})
    \\[1mm]& <\Big(\frac{\rho_n}{c_n} -2C\epsilon + V (x_{\lambda,\epsilon})\Big) w (x_{\lambda,\epsilon})+w  (x_{\lambda,\epsilon})  \int_{H_{\lambda}} \frac{ 1_{B_1(x_{\lambda,\epsilon})}(y )+1_{B_1(x_{\lambda,\epsilon})}(y^{\lambda}) }{\left|x_{\lambda,\epsilon}-y^{\lambda}\right|^{n}  } dy,
\end{align*}
which implies that
$$\int_{H_{\lambda}} \frac{ 1_{B_1(x_{\lambda,\epsilon})}(y )+1_{B_1(x_{\lambda,\epsilon})}(y^{\lambda}) }{\left|x_{\lambda,\epsilon}-y^{\lambda}\right|^{n}  } dy<-V (x_{\lambda,\epsilon})+\frac{\rho_n}{c_n} +2C\epsilon \leq C'.$$
This contradicts  (\ref{lim-oo-1}) for $\epsilon>0$ small enough.
Therefore,  the negative minimal of $w$   can't be achieved  in  $\Omega\cap    \bP_{\lambda,\epsilon} $.   \hfill$\Box$\medskip

\subsection{ The moving plane method }

\noindent {\bf Proof of Theorem \ref{teo 1sym}.}   Let $v$ be a positive solution of (\ref{eq 1.1sym}) 
and $v^\#$ be  the Kelvin transform
of $v$ centered at $\tilde x$, defined in (\ref{def k-0}).   Then  it follows by Proposition \ref{cr 3.0} with $T=\rho_n$ that
 \begin{equation*}
\left\{ \arraycolsep=1pt
\begin{array}{lll}
\cL_\Delta v^{\#_{0,1}}- \rho_nv^{\#_{0,1}} = \frac4n v^{\#_{0,1}} (\ln v^{\#_{0,1}}) \quad \  &{\rm in}\ \   \R^n\setminus \{\tilde x\},\\[2mm]
 \phantom{\cL_\Delta  - \rho_nv^{\#_{0,1}}  }
  v^{\#_{0,1}} \geq 0\quad \ &{\rm{in}}  \ \ \R^n.
\end{array}
\right.
\end{equation*}
For simplicity, we set $\tilde x=0$, $u_1= v^{\#_{0,1}}$  and  in the following, we  employ the moving plane method to show (\ref{con h}).

We rewrite the equation (\ref{eq 1.1sym}) as follows.
\begin{equation}\label{eq 1.1sym-re}
\left\{ \arraycolsep=1pt
\begin{array}{lll}
\cL_\Delta u_1  - \rho_n u_1+f_1(u_1)= f_2(u_1)\quad \  &{\rm in}\ \   \R^n\setminus \{0\},\\[2mm]
 \phantom{\cL_\Delta  +(\rho- \rho_n)u_l  }
u_1> 0\quad \ &{\rm{in}}  \ \ \R^n\setminus \{0\} ,
\end{array}
\right.
\end{equation}
 where $f_1, \ f_2$ are defined by
$$f_1(t)=\left\{ \arraycolsep=1pt
\begin{array}{lll}
 - \frac4n t\ln t \quad \  &{\rm for}\quad  t\in [0,\frac{1}{e}],\\[2mm]
 \phantom{  }
 \frac{4}{n e}  \quad \ &{\rm{for}}  \quad t\geq \frac1e
\end{array}
\right. $$
and
$$f_2(t)=\left\{ \arraycolsep=1pt
\begin{array}{lll}
0 \quad \  &{\rm for}\quad  t\in[0,\frac1e],\\[2mm]
 \phantom{  }
\frac4n\big( t\ln t +   \frac1e\big)  \quad \ &{\rm{for}}  \quad t\geq \frac1e.
\end{array}
\right. $$
Let
$$\cO_t=\big\{x\in\R^d:\,  u_1(x)>t \big\}\quad {\rm for}\ \, t>0, $$
then
$$f_1(u_1)\equiv  \frac{4}{n e} \quad {\rm in}\ \,   \cO_{\frac1{e}}\quad{\rm and}\quad  f_2(u_1)\equiv0\quad {\rm in}\ \, \R^n\setminus \cO_{\frac1{e}}. $$
Form (\ref{def k-0}), we note that
$$
 \lim_{|x-\tilde x|\to+\infty} u_1(x)|x|^n= v^\#(0)\in(0,+\infty).
$$


For $\lambda< 0$,  denote
 $$
u_{1,\lambda}(x) =u_1(x^\lambda),\qquad  w_{\lambda}(x)=u_{1,\lambda}(x)-u_1(x),
$$
then $w_{\lambda}$ verifies the following equation
\begin{equation}\label{eq 1.1sym-re-}
\left\{ \arraycolsep=1pt
\begin{array}{lll}
\cL_\Delta w_{\lambda} -\rho_n w_{\lambda}+ V_{\lambda} w_{\lambda}=W_{\lambda} w_{\lambda}    \quad \  &{\rm in}\ \,   H_\lambda,
\\[2.5mm]
 \phantom{ \cL_{\Delta}  n) V_{l,\lambda}-w_{\lambda}\, }
  w_{\lambda} ( x^\lambda )=-w_{\lambda} (x) &   {\rm in }\ \,  H_\lambda,
\end{array}
\right.
\end{equation}
where
\begin{equation} \label{V-lambda}
V_{\lambda} (x)=\left\{ \arraycolsep=1pt
\begin{array}{lll}\frac{f_1(u_{1,\lambda}(x))-f_1(u_1(x))}{u_{1,\lambda}(x)-u_1(x)}\quad &{\rm if}\ \, u_{1,\lambda}(x)\not=u_1(x),
\\[2mm]
f_1'(u_1(x)) \quad  & {\rm if}\ \, u_{1,\lambda}(x) =u_1(x)
\end{array}
\right.
\end{equation}
and
\begin{equation} \label{W-lambda}
W_{\lambda} (x)=\left\{ \arraycolsep=1pt
\begin{array}{lll}\frac{f_2(u_{1,\lambda}(x))-f_2(u_1(x))}{u_{1,\lambda}(x)-u_1(x)}\quad &{\rm if}\ \, u_{1,\lambda}(x)\not=u_1(x),
\\[2mm]
f_2'(u_1(x)) \quad  & {\rm if}\ \, u_{1,\lambda}(x) =u_1(x),
\end{array}
\right.
\end{equation}
 which are nonnegative and  continuous in $H_\lambda$.

Observe that  $w_{\lambda}$ vanishes at infinity and on the boundary of $H_\lambda$.
 If
 $$\inf_{z\in H_\lambda}w_{\lambda}(z)<0,$$
  then   there exists  $x_{0,\lambda}\in H_\lambda$ such that
  $$w_{\lambda}(x_{0,\lambda}) =\inf_{z\in H_\lambda}w_{\lambda}(z). $$
Denote
$$\Omega_\lambda=\big\{z\in H_\lambda:\, w_{\lambda}(z)<0\big\}. $$
  \smallskip

{\it Step 1:  Start moving the plane $T_{\lambda}$ along the $x_{1}$-axis from near $-\infty$ to the right.}  We need to show
that $w_{\lambda}>0$  in $H_\lambda$.

Since $\displaystyle \lim_{|x|\to+\infty}u_1(x)=0$ and $u_{1,\lambda}(x)<u_1(x)$ in $\Omega_\lambda$, then there exists $\bar \lambda<0$ such that
for any $\lambda<\bar \lambda$, $u_1(x)<\frac1e$. Then we have that
$$f_2(u_1(x))=f_2(u_{1,\lambda})=0\quad {\rm for}\ \ x\in\Omega_\lambda  $$
and
 $w_{\lambda}$ verifies the equation
\begin{equation}\label{eq 1.1sym-re}
\left\{ \arraycolsep=1pt
\begin{array}{lll}
\cL_\Delta w_{\lambda} -\rho_n w_{\lambda}+ V_{\lambda} w_{\lambda}=0    \quad \  &{\rm in}\ \,   \Omega_\lambda,
\\[2.5mm]
 \phantom{ \cL_{\Delta}   -\rho_n)  -w_{\lambda}-\! }
w_{\lambda}(x)\geq 0\quad &{\rm in}\ \,  H_\lambda\setminus \Omega_\lambda,\\[2.5mm]
 \phantom{ \cL_{\Delta} \rho_n) V_{1\lambda}-w_{\lambda}}
  w_{\lambda}( x^\lambda )=-w(x)\qquad  &    {\rm in }\ \,  H_\lambda,
\end{array}
\right.
\end{equation}
where $V_{\lambda}$ is defined as in (\ref{V-lambda}).

Let
$$
w_{\lambda,-}(x)=\left\{ \arraycolsep=1pt
\begin{array}{lll}
\min\{0,w_{\lambda}(x)\}\quad &{\rm for}\ \, x\in H_\lambda,\\[2.5mm]
\max\{0,w_{\lambda}(x)\} \quad &{\rm for}\ \, x\in \R^n\setminus H_\lambda,
\end{array}
\right.
$$
which belongs to $\cD^0(\R^n)\cap L_0(\R^n)$. Multiply $w_{\lambda,-}$ and integrate over $\R^n$ to obtain that
  \begin{align}\label{inequa-1}
  \int_{\R^n}  (\cL_{\Delta} w_{\lambda}  )w_{\lambda,-} dx =\int_{\R^n}  \big(\rho_n- V_{\lambda} \big) |w_{\lambda,-}|^2 dx.
  \end{align}
Furthermore, we want to employ the following a sharp form of Pitt's inequality
 \cite[Theorem 1]{Beckner} that
  $$\int_{\R^n} (\ln |x| )|f(x)|^2dx+\int_{\R^n} (\ln |\xi| )|\hat{f}(\xi)|^2dx\geq D_n\int_{\R^n} |f(x)|^2dx,  $$
  where $D_n=\psi(\frac{n}{4})-\ln \pi$.

  Moreover, we obtain that  for $x\in {\rm supp}(w_{\lambda,-})$,
  \begin{align*}
 \frac1{c_n} \Big(\cL_\Delta w_{\lambda,-}(x) - \cL_\Delta w_{\lambda}(x)\Big)&= \int_{\R^n} \frac{w_{\lambda}(z)-w_{\lambda,-}(z) }{|x-z|^n}dz
  \\[1mm]&=  \int_{H_\lambda \cap ({\rm supp}(w_{\lambda,-}))^c}  \frac{w_{\lambda}(z) }{|x-z|^n}dz+ \int_{H_\lambda^c \cap ({\rm supp}(w_{\lambda,-}))}  \frac{w_{\lambda}(z) }{|x-z|^n}dz
   \\[1mm]&=  \int_{H_\lambda \cap ({\rm supp}(w_{\lambda,-}))^c} w_{\lambda}(z)\Big( \frac{1}{|x-z|^n} - \frac{1}{|x-z_\lambda |^n}  \Big)dz
     \\[1mm]& \geq0,
  \end{align*}
where $w_{\lambda}(z)>0$ for $z\in H_\lambda \cap ({\rm supp}(w_{\lambda,-}))^c$.  That is
    \begin{align*}
 \cL_\Delta w_{\lambda,-}(x) \geq  \cL_\Delta w_{\lambda}(x),\quad\forall\,  x\in {\rm supp}(w_{\lambda,-}).
  \end{align*}
   Thus,  together with (\ref{inequa-1}) we obtain that
  \begin{align*}
  \int_{\R^n}  \big(\rho_n- V_{\lambda} \big) |w_{\lambda,-}|^2 dx&=\int_{\R^n\cap( {\rm supp}( w_{\lambda,-} )) }  (\cL_{\Delta} w_{\lambda}  )w_{\lambda,-} dx
  \\&\geq \int_{\R^n}  (\cL_{\Delta} w_{\lambda,-}  )w_{\lambda,-}dx
    \\& =\int_{\R^n}  (2\ln|\xi|) |\hat{w}_{\lambda,-}|^2 dx
     \\&\geq 2D_n\int_{\R^n} |w_{\lambda,-}(x)|^2dx-2\int_{\R^n} |w_{\lambda,-}(x)|^2 \ln|x|dx,
\end{align*}
that is
  \begin{align}\label{ineq-mo-1}
  \int_{H_\lambda}  \big(2\ln |x|- V_{\lambda} \big) |w_{\lambda,-} |^2 dx \geq  (2D_n-\rho_n)\int_{H_\lambda} w_{\lambda,-}(x)^2dx ,
\end{align}
where
 \begin{align*}
 V_{\lambda}(x)=   \frac{f_1(u_{1, \lambda}(x))-f_1(u_{1}(x))}{u_{1,\lambda}(x)-u_{1}(x)}
\geq  -\frac4n\big(\ln (u_1(x))+1\big)
\geq 4 \ln |x|-\frac4n -\ln (2c),
\end{align*}
by the assumption that  $\displaystyle  \lim_{|x| \to+\infty} u_{1}(x)|x|^n=c>0$.

Thus, (\ref{ineq-mo-1}) reduces to  that for $\lambda\ll-1$ and $C>0$ independent of $\lambda$
  \begin{align}\label{ineq-mo-2}
  -2\int_{H_\lambda} (\ln |x|)  w_{\lambda,-} ^2 dx \geq  (2D_n-\rho_n+C)\int_{H_\lambda} w_{\lambda,-}^2dx,
\end{align}
which leads to $w_{\lambda,-}\equiv 0$ in $H_\lambda$ for $\lambda\leq \tilde \lambda$,
where  $\tilde \lambda<0$ is negative enough.

So $w_{\lambda}\geq  0$ in $H_\lambda$ and it follows  by Lemma \ref{strong} that $w_{\lambda}>  0$ in $H_\lambda$,
since $w_{\lambda}(0^\lambda)= u_1(0)-u_1(0^\lambda)>0$,
where $u_1(0)>0$ and $u_1(0^\lambda)\to0$ as $\lambda\to-\infty$.
 \smallskip

 {\it Step 2: Keep moving the plane to the limiting position $T_0$.}
Let
$$
\lambda_{0}=\sup \big\{\lambda<0\!:  \,  w_{1,\mu}(x) \geq 0,\,  \forall x \in H_{\mu},\, \mu \leq \lambda\big\}.
$$
In this part, we show that
\begin{equation}\label{e 3.10}
\lambda_{0}=0.
\end{equation}

Otherwise,  supposing that $\lambda_{0}<0$, the plane $T_{\lambda_{0}}$ could still be moved further to the right.
By the continuity and strong maximum principle, we have that
$$
{\rm either}\qquad  w_{\lambda_0}(x)\equiv 0,\ \,\forall\, x \in  H_{\lambda_0} \qquad{\rm or}\qquad w_{\lambda_0}(x)> 0,\ \,\forall\, x \in  H_{\lambda_0}.
$$

If $\lambda_0<0$ and $$w_{\lambda_0}(x)\equiv 0,\ \,\forall\, x \in  H_{\lambda_0}, $$
then $u_1(0)= u_1(0^{\lambda_0})<+\infty$,
which implies (\ref{con h}).  We are done.


So in the latter proof, we can alway assume that  $$w_{\lambda_0}(x)> 0,\ \,\forall\, x \in  H_{\lambda_0}, $$
if $\lambda_0<0$.  In this case, we continue to move the  planes to the right.

Now take  $\lambda_\epsilon=\lambda_0+\epsilon$ with $\epsilon\in(0,-\frac{1}{4}\lambda_0)$.  Let $\epsilon_0\in(0,\frac{1}{8}\min\{1,-\lambda_0\})$ and   for $\epsilon\in(0,\epsilon_0)$ there exists $R_\epsilon>4(R_1-\lambda_0)$ such that
$R_\epsilon\to+\infty$ as $\epsilon\to0^+$ and
$$
w_{\lambda_\epsilon}(x)> 0,\ \,\forall\, x \in  H_{\lambda_0-\epsilon} \cap B_{R_\epsilon}(0).
$$
where $R_1>1$ such that $$\cO_{\frac1{e}}\subset B_{ R_1}(0).  $$
That means
$$\Omega_{\lambda_\epsilon}\subset   \big(H_{\lambda_\epsilon}\setminus H_{\lambda_0-\epsilon}\big)  \cup \big(H_{\lambda_\epsilon} \setminus B_{R_\epsilon}(0)\big),  $$
where
$\Omega_{\lambda_\epsilon}=\{ x\in H_{\lambda_\epsilon}: w_{1,\lambda_\epsilon}(x)<0\}$.
Since $\displaystyle \lim_{|x'|\to+\infty}u_1(x_1,x')=0$, so is $w_{\lambda}$.  If $\Omega_{\lambda_\epsilon}\not=\emptyset$, then
the minimal of $w_{1,\lambda_\epsilon}$ in $H_{\lambda_\epsilon}$ can be achieved by some point $x_{0,\epsilon} \in H_{\lambda_\epsilon}$.

From Lemma \ref{Narrow A}, the minimal can't be achieved in $H_{\lambda_\epsilon}\setminus  H_{\lambda_0-\epsilon}  $.  Now we assume that the minimal is achieved in $H_{\lambda_0-\epsilon} \setminus   B_{R_\epsilon}(0) $. So there exists $x_{0,\epsilon}\in H_{\lambda_0-\epsilon} \setminus  B_{R_\epsilon}(0) $ such that
$$w_{\lambda_\epsilon} (x_{0,\epsilon}) =\inf_{x\in \Omega_{\lambda_\epsilon}} w_{\lambda_\epsilon} (x)<0. $$

 Remark that for $\epsilon >0$ small enough
 $$u_{1,{\lambda_\epsilon}}<\frac1{e}\quad \text{ in $H_{\lambda_\epsilon}\setminus B_{R_\epsilon}(0)$.}$$
Since
$$u_{1,\lambda_\epsilon}(x_{0,\epsilon})<  u_1(x_{0,\epsilon})  <\frac1e,$$
 then
 $$f_2(u_1(x_{0,\epsilon}))=f_2(u_{1,\lambda_\epsilon}(x_{0,\epsilon}))=0. $$

 Note that
\begin{equation}\label{eq 1.1sym-re}
\left\{ \arraycolsep=1pt
\begin{array}{lll}
\cL_{\Delta} w_{\lambda_\epsilon} - \rho_n  w_{\lambda_\epsilon} +   V_{\lambda_\epsilon} w_{\lambda_\epsilon}  =0 \quad \  &{\rm in}\ \,   \Omega_{\lambda_\epsilon},
\\[2.5mm]
 \phantom{ \cL_{\Delta,\R^n} - }
  w_{\lambda_\epsilon} ( x^{\lambda_\epsilon}  )=-w_{\lambda_\epsilon} (x) &   {\rm in }\ \,   H_{\lambda_\epsilon},\\[2.5mm]
   \phantom{ \cL_{\Delta,\R^n}  +\    }
  w_{\lambda_\epsilon} (x)\geq 0 &   {\rm in }\ \,    H_{\lambda_\epsilon}\setminus \Omega_{\lambda_\epsilon},
\end{array}
\right.
\end{equation}
where
$$V_{\lambda_\epsilon} (x)=\left\{ \arraycolsep=1pt
\begin{array}{lll}
\frac4n \frac{f_1(u_{1,\lambda_\epsilon}(x))-f_1(u_1(x))}{u_{1,\lambda_\epsilon}(x)-u_1(x)}\quad\ \ &{\rm if}\ \, u_{1,\lambda}(x)\not=u_1(x),
\\[2mm]
\frac4n f_1'(u_1(x)) \quad  & {\rm if}\ \, u_{1,\lambda}(x) =u_1(x),
\end{array}
\right.
$$
which is nonnegative and  continuous in $(0, +\infty)$.

Since $|x_{0,\epsilon}|>R_\epsilon\to+\infty$ as $\epsilon\to0^+$, then
 $$0<u_{1,\lambda_\epsilon}(x_{0,\epsilon})< u_1(x_{0,\epsilon})\to0\quad \text{ as $x_{0,\epsilon}\to+\infty$}$$
 and
 \begin{align*}
 V_{\lambda_\epsilon}(x_{0,\epsilon}) &=\frac4n\frac{f_1(u_{1,\lambda_\epsilon}(x_{0,\epsilon}))-f_1(u_1(x_{0,\epsilon}))}{u_{1,\lambda_\epsilon}(x_{0,\epsilon})-u_1(x_{0,\epsilon})}
 \\[1mm]&\geq\frac4n f_1'(u_1(x_{0,\epsilon}))
 =-\frac4n\big(\ln u_1(x_{0,\epsilon})+1\big)\to+\infty\quad{\rm as} \ \epsilon\to0^+.
  \end{align*}


Observe that
 \begin{align}\label{aa-00-1}
\cL_\Delta  w_{\lambda_\epsilon}-\rho_n w_{\lambda_\epsilon} +  V_{{\lambda_\epsilon}}  w_{\lambda_\epsilon}   =  \cL_{1,\epsilon}w_{\lambda_\epsilon}(x)+\cL_{2,\epsilon}w_{\lambda_\epsilon}(x), \quad\forall\,  x\in\Omega_{\lambda_\epsilon},
 \end{align}
where
$$\cL_{1,\epsilon}w_{\lambda_\epsilon}(x)= c_n \int_{\R^n} \frac{w_{\lambda_\epsilon} (x )-w_{\lambda_\epsilon} (y)}{ |x -y|^{n}} 1_{B_1(x)}(y)  d y$$
and
$$\cL_{2,\epsilon}w_{\lambda_\epsilon}(x)=    V_{\lambda_\epsilon}w_{\lambda_\epsilon}(x)- c_n \int_{\R^n}\frac{w_{\lambda_\epsilon} (y)}{|x -y|^{n}}1_{B_1^c(x)}(y) dy.
$$

{\bf Claim 1: }  There exists $ \epsilon_0>0$ such that for $ \epsilon \in(0, \epsilon_0)$
$$  \cL_{2,\epsilon}(x_{0,\epsilon})\leq0. $$
 In fact, we see that
 \begin{align*}
\frac{1}{c_n} \cL_{2,\epsilon}w_{\lambda_\epsilon}(x_{0,\epsilon})&\geq  \frac{1}{c_n}  V_{\lambda_\epsilon}(x_{0,\epsilon})   w_{\lambda_\epsilon}  (x_{0,\epsilon}) -\Big( \int_{ H_{\lambda_\epsilon} } \frac{ w_{\lambda_\epsilon} (y)1_{B_1^c(x_{0,\epsilon})}(y) }{\left|x_{0,\epsilon}-y\right|^{n}}  d y+\int_{ \tilde  H_{\lambda_\epsilon} } \frac{ w_{\lambda_\epsilon} (y)1_{B_1^c(x_{0,\epsilon})}(y)}{\left|x_{0,\epsilon}-y\right|^{n}}  d y \Big)
 \\[1mm]&=\frac{1}{c_n}  V_{\lambda_\epsilon}(x_{0,\epsilon})     w_{\lambda_\epsilon} (x_{0,\epsilon}) -\Big( \int_{ H_{\lambda_\epsilon} } \frac{ w_{\lambda_\epsilon} (y)1_{B_1^c(x_{0,\epsilon})}(y) }{\left|x_{0,\epsilon}-y\right|^{n}}  d y-\int_{  H_{\lambda_\epsilon} } \frac{ w_{\lambda_\epsilon} (y)1_{B_1^c(x_{0,\epsilon})}(y^\lambda) }{\left|x_{0,\epsilon}-y^{\lambda_\epsilon} \right|^{n}} d y \Big)
 \\[1mm]&=\frac{1}{c_n}  V_{\lambda_\epsilon}(x_{0,\epsilon})     w_{\lambda_\epsilon} (x_{0,\epsilon})- \int_{ H_{\lambda_\epsilon} } w_{\lambda_\epsilon}(y) 1_{B_1^c(x_{0,\epsilon})}(y)   \Big( \frac{1}{\left|x_{0,\epsilon}-y\right|^{n}}- \frac{1}{\left|x_{0,\epsilon}-y^{\lambda_\epsilon} \right|^{n}}\Big) dy
 \\[1mm]&\qquad \qquad\  -\int_{ H_{\lambda_\epsilon} } w_{\lambda_\epsilon}  (y)  \frac{1_{B_1^c(x_{0,\epsilon})}(y)-1_{B_1^c(x_{0,\epsilon})}(y^{\lambda_\epsilon})   }{\left|x_{0,\epsilon}-y^{\lambda_\epsilon}\right|^{n}}   dy,
\end{align*}
where  $B_1^c(x_{0,\epsilon})=\R^n\setminus B_1(x_{0,\epsilon})$, $x_{0,\epsilon}=((x_{0,\epsilon})_1, x_{0,\epsilon}')\in (0,\lambda_\epsilon)\times \R^{n-1}$,
\begin{align*}
 \Big|\int_{ H_{\lambda_\epsilon} } w_{\lambda_\epsilon}  (y)  \frac{1_{B_1^c(x_{0,\epsilon})}(y)-1_{B_1^c(x_{0,\epsilon})}(y^{\lambda_\epsilon})   }{\left|x_{0,\epsilon}-y^{\lambda_\epsilon}\right|^{n}}   dy\Big|&\leq  -w_{\lambda_\epsilon}(x_{0,\epsilon})   \int_{ \big( B_1(x_{0,\epsilon})\setminus B_1(x_{0,\epsilon}^{\lambda_\epsilon})\big) \cap H_{\lambda_\epsilon}  }  \frac{ 1  }{\left|x_{0,\epsilon}-y^{\lambda_\epsilon}\right|^{n}}   dy
 \\[1mm]&\leq -w_{\lambda_\epsilon}(x_{0,\epsilon})  |B_1(x_{0,\epsilon})\setminus B_1(x_{0,\epsilon}^{\lambda_\epsilon})|
 \\[1mm]& \leq -w_{\lambda_\epsilon}(x_{0,\epsilon}) |B_1(0)|
 \end{align*}
and
 \begin{align*}
 &\quad\ -\int_{ H_{\lambda_\epsilon} } w_{\lambda_\epsilon}(y) 1_{B_1^c(x_{0,\epsilon})}(y)   \Big( \frac{1}{\left|x_{0,\epsilon}-y\right|^{n}}- \frac{1}{\left|x_{0,\epsilon}-y^{\lambda_\epsilon} \right|^{n}}\Big) dy
\\[1mm] &\leq -w_{\lambda_\epsilon}(x_{0,\epsilon})\int_{ H_{\lambda_\epsilon} }1_{B_1^c(x_{0,\epsilon})}(y)   \Big( \frac{1}{\left|x_{0,\epsilon}-y\right|^{n}}- \frac{1}{\left|x_{0,\epsilon}-y^{\lambda_\epsilon} \right|^{n}}\Big) dy
\\[1mm] &\leq -w_{\lambda_\epsilon}(x_{0,\epsilon}) \int_{ H_{\lambda_\epsilon} }\frac{1_{B_1^c(x_{0,\epsilon})}(y) }{\left|x_{0,\epsilon}-y\right|^{n-1}}   \Big(\frac{1}{\left|x_{0,\epsilon}-y\right|} - \frac{1}{\left|x_{0,\epsilon}-y^{\lambda_\epsilon} \right| }\Big) dy
\\[1mm] &\leq -w_{\lambda_\epsilon}(x_{0,\epsilon}) \int_{ H_{\lambda_\epsilon} }\frac{1_{B_1^c(x_{0,\epsilon})}(y) }{\left|x_{0,\epsilon}-y\right|^{n+1}}   \Big( |x_{0,\epsilon}-y^{\lambda_\epsilon} |- |x_{0,\epsilon}-y |  \Big)  dy
\\[1mm] &\leq -w_{\lambda_\epsilon}(x_{0,\epsilon}) \int_{ H_{\lambda_\epsilon} }\frac{1_{B_1^c(x_{0,\epsilon})}(y) }{\left|x_{0,\epsilon}-y\right|^{n+1}}   \Big( |x_{0,\epsilon}^{\lambda_\epsilon} -y |- |x_{0,\epsilon}-y |  \Big)  dy
\\[1mm] &= -w_{\lambda_\epsilon}(x_{0,\epsilon}) \int_{ H_{\lambda_\epsilon} }\frac{1_{B_1^c(x_{0,\epsilon})}(y) }{\left|x_{0,\epsilon}-y\right|^{n+1}}      dy\,  |x_{0,\epsilon} -x_{0,\epsilon}^{\lambda_\epsilon}|
\\[1mm] &= -2(\lambda_\epsilon-(x_{0,\epsilon})_1) w_{\lambda_\epsilon}(x_{0,\epsilon})\Big( \int_{B_1(0)^c}\frac{1 }{ |z|^{n+1}}      dz\Big)\,
\\[1mm] &\leq -2(\lambda_0+1) w_{\lambda_\epsilon}(x_{0,\epsilon})\Big( \int_{B_1(0)^c}\frac{1 }{ |z|^{n+1}}      dz\Big).
\end{align*}

Thus, we conclude that
 \begin{align*}
 \cL_{2,\epsilon}w_{\lambda_\epsilon}(x_{0,\epsilon})&\leq \Big(  V_{\lambda_\epsilon}(x_{0,\epsilon})    -2c_n (\lambda_0+1)   \int_{B_1(0)^c}\frac{1 }{ |z|^{n+1}}      dz-|B_1(0)|\Big)w_{\lambda_\epsilon}  (x_{0,\epsilon})\leq 0
\end{align*}
 if
$$    V_{\lambda_\epsilon}(x_{0,\epsilon}) \geq  2c_n (\lambda_0+1)    \Big( \int_{B_1(0)^c}\frac{1 }{ |z|^{n+1}}      dz\Big)+c_n|B_1(0)|, $$
 which could be guaranteed for $\epsilon>0$ small enough.

As it is shown in (\ref{pp-1-1}) that
\begin{align}\label{pp-1+1}
   \cL_{1,\epsilon}w_{\lambda_\epsilon}(x_{\lambda,\epsilon})  <         w_{\lambda_\epsilon}(x_{0,\epsilon})\int_{H_{\lambda}} \frac{ 1_{B_1(x_{0,\epsilon})}(y )+1_{B_1(x_{0,\epsilon})}(y^{\lambda_\epsilon}) }{\left|x_{0,\epsilon}-y^{\lambda_\epsilon}\right|^{n}  } dy.
\end{align}

 As a consequence,  together with {\bf Claim 1} and (\ref{pp-1+1}), we obtain that
 \begin{align}\label{baa-00-2}
0=\frac1{ c_{n}}\Big(\cL_\Delta  w_{\lambda_\epsilon}   +c_n V_{{\lambda_\epsilon}}  w_{\lambda_\epsilon} \Big)(x_{0,\epsilon}) <w_{\lambda_\epsilon} (x_{0,\lambda})  \int_{ H_{\lambda_\epsilon} } \frac{ 1_{B_1(x_{0,\lambda})}(y )+1_{B_1(x_{0,\lambda})}(y^{\lambda_\epsilon} ) }{\left|x_{0,{\lambda_\epsilon} }-y^{\lambda_\epsilon} \right|^{n}  } dy,
 \end{align}
which deduces  that
\begin{align}\label{bound-3.1}
  \int_{ H_{\lambda_\epsilon} } \frac{ 1_{B_1(x_{0,\epsilon})}(y )+1_{B_1(x_{0,\epsilon})}(y^{\lambda_\epsilon} ) }{\left|x_{0,\epsilon}-y^{\lambda_\epsilon} \right|^{n}  } dy<0,
\end{align}
which is impossible.

Therefore,  $\lambda_0=0$ and $w_{\lambda_0}\geq0$ in $ H_{0}$.
 That is
 $$u_1(x_1,x')\leq u_1(-x_1,x')\quad \text{ for $(x_1,x')\in  H_0 $}.$$ \smallskip

 For $\lambda \geq 0$, we move the plane from near $x_{1}=+\infty$ to the left, and by a similar argument, we can derive that
$$
u_1\left(-x_{1}, x'\right) \geq u_1\left(x_{1}, x'\right), \quad x_{1}>0,\ x'\in\R^{n-1}.
$$
Then
$$
u_1\left(-x_{1}, x'\right) = u_1\left(x_{1}, x'\right).
$$
Now we derive that $u_1(x)$ is symmetric about the plane $T_{0}$.  In fact, we can do the moving planes
in any direction to the position $T_{0}$,  thus we have that
\begin{align}\label{bound-3.1---}
u_1(x)=u_1( y)\quad {\rm for}\ x,y\in\R^{n}, \,  |x|=|y|.
\end{align}

If $u_1$ is bounded at  $\tilde x$, then the radial symmetry and nondecreasing monotonicity lead
to
$$u_1(\tilde x)=\lim_{|x-\tilde x|\to 0} u_1(x)\in(0,+\infty)  $$
and
$$\lim_{|x|\to+\infty} u_1(x)|x|^n=v(\tilde x)$$
which implies (\ref{con h})
and we are done.

If $u_1$ is unbounded at  $x$, by the symmetry,  Therefore, $v$ blows up at  $\tilde x$, which implies that
 $$\lim_{|x|\to+\infty} v(x)|x|^n=+\infty. $$
  Now for any  two different points $x,\, \tilde x$  in $\R^n$, let
$\bar x=\frac{\tilde x+x}{2}$, which is the midpoint of the line segment $\overline{x,\tilde x}$, then from our method of moving planes for $\bar u_1$ to the position $\bar x$, we obtain that
$v(x)=v(\tilde x)$, since $x,\tilde x$ are arbitrary  two points in $\R^n$ , then $v$ must be a positive constant and it is impossible for $v\in L^1_0(\R^n)$. \smallskip

As a consequence,  (\ref{con h}) holds ture, i.e.
 $$\lim_{|x|\to+\infty} v(x)|x|^n\in(0,+\infty). $$

{\it  Radial symmetry  and decreasing monotonicity of $v$.} Under the assumption the decay at infinity (\ref{con h}),   the radial symmetry about some point and decreasing monotonicity follow the above method of the moving planes again direct for the solution $v$ of  the problem (\ref{eq 1.1}). \hfill$\Box$ \medskip

\begin{remark}\label{sect 4-rm1}
$(i)$ From the asymptotic behavior (\ref{con h}), we have that $u\in L^q(\R^n)$ for  any positive solution $u$ of (\ref{eq 1.1}) and any $q>1$.

$(ii)$ Our procedure of moving plane depends on  nonlinearity $ \frac4n u \ln u$, which is essential   for the starting the moving planes from the infinity.

\end{remark}

  \setcounter{equation}{0}
\section{Classification}

In this section,  we   focus on the uniqueness of the solution of  \eqref{eq 1.1}.
      \begin{theorem}\label{teo 1}
       Let $u\in L^1_0(\R^n)\cap \cD^0(\R^n)$   be a  positive solution $u$ of \eqref{eq 1.1},
then there exist   $t>0$ and $\tilde x\in\R^n$ such that
   \begin{equation}\label{sol form-1}
   u(x)= u_t(x-\tilde x)\quad\ {\rm for}\ \, x\in\R^n,
    \end{equation}
    where $u_t$ is defined in (\ref{sol form-00}).
    \end{theorem}

\subsection{Uniqueness under some restriction}
Note that  $u_t$ with $t>0$,  defined in (\ref{sol form-0}),  is the solution of \eqref{eq 1.1} with the invariant $L^2(\R^n)$ norm,
  that is,
   \begin{equation}\label{Norm 1}
\Lambda_n:=   \|u_t\|_{L^2(\R^n)} =\beta_n \sqrt{ \int_{\R^n}\frac{1}{(1+|x|^2)^n} dx} =\beta_n \pi^{\frac n4} \sqrt{ \frac{\Gamma(\frac n2)}{\Gamma(n)}}.
   \end{equation}

\begin{proposition}\label{teo 1-res}
       Let  $u_t$ be defined in (\ref{sol form-0}) and  $\Lambda_n$ be the constant defined in (\ref{Norm 1}).

  If   $u$  is a  positive solution  of \eqref{eq 1.1}
   such that
       \begin{equation}\label{res-1}
       \|u\|_{L^2(\R^n)}= {\Lambda_n},
        \end{equation}
 then there exist  $t>0$ and $\tilde x\in\R^n$ such that
 $$
   u(x)= u_t(x-\tilde x)\quad\ {\rm for}\ \, x\in\R^n. $$
  \end{proposition}
We will adapt the idea of \cite{CLO}  to show Proposition \ref{teo 1-res}. So we provide the following auxiliary lemmas.
 Let
 $$u_1(x)=\beta_n \frac{1}{(1+|x|^2)^{\frac n2}}, $$
 which is the standard solution of \eqref{eq 1.1}  centered at the origin.

\begin{lemma} \label{lm 5.1-1}
Let $u$  be a  positive solution  of \eqref{eq 1.1}, $x_0\in\R^n$ and
$$s^n=\frac{u_\infty}{u(x_0)}. $$
Then
$$ u(sx+x_0)=|x|^{-n} u\big(\frac{sx}{|x|^2}+x_0\big), \quad \forall x\in\R^n. $$
\end{lemma}
\noindent{\bf Proof. } Let $x_0=0$,  $e$ be a unit vector in $\R^n$ and
$$w(x)=|x|^{-n} u(\frac{x}{|x|^2}-se), $$
then
\begin{align*}
w(0)=\lim_{|x|\to0^+} w(x)&=\lim_{|x|\to0^+} |x|^{-n} u(\frac{x}{|x|^2}-se)
\\[1mm]&= \lim_{|x|\to+\infty} |z|^{n} u(z-se)=u_\infty=s^{n}u(0)=w(\frac1s e)
\end{align*}
and $w$ is symmetric about $\frac{1}{2s} e$.
For $h>0$,  let
$$z_1=\frac{\frac12-h}{s}e\quad {\rm and}\quad z_2=\frac{\frac12+h}{s}e, $$
then $z_1+z_2=\frac{1}{s} e$,
$$s\frac{\frac12+h}{\frac12-h}e =\frac{z_1}{|z_1|^2}-se,\qquad s\frac{\frac12-h}{\frac12+h}e=\frac{z_2}{|z_2|^2}-se$$
and
\begin{align*}
 s^{n}|\frac1{2s}-h|^{-n} u(s\frac{\frac12+h}{\frac12-h}e)  =  w(z_1)&= w(z_2)= s^{n}|\frac1{2s}-h|^{-n} u(s\frac{\frac12-h}{\frac12+h}e),
\end{align*}
which, letting $t=\frac{\frac12-h}{\frac12+h}$,  implies that
$$u(ste)=t^{-n} u(st^{-1}e). $$
We end the proof from a translation if $x_0\not=0$.
  \hfill$\Box$\smallskip

   \begin{lemma} \label{lm 5.2-1}
 Let $w$ be a positive solution  of \eqref{eq 1.1}  such that $w$   verifies (\ref{res-1}) and is  radial symmetric about $x_0$,  decreasing in the radial direction $|x-x_0|$,
 then
$$ w(x_0) w_\infty=\beta_n^2. $$
\end{lemma}
\noindent{\bf Proof.}
By a translation and a rescaling in Theorem \ref{cr 3.0},  we can assume that $x_0=0$  and $w_\infty=\beta_n$.  Next we need to show $w(0)=\beta_n$.

{\it  We first deal with
the case that $w(0)>\beta_n$. } To this end,  we denote
$$w_0(\iota)=w(\iota e)\quad{\rm and}\quad  u_0(\iota)=u_1(\iota e)\ \ {\rm for}\ \, \iota\in\R,$$
 where $e$ is a unit vector in $\R^n$.  By the radial symmetries,  the functions $w$ and $u_0$ are independent of $e$.

When $w(0)>\beta_n$,   we    assume  that  there exists $\iota_0>0$ such that
\begin{equation}\label{a--}
w_0(\iota)>u_0(\iota) \quad {\rm for}\ |\iota|<\iota_0\ \ {\rm and}\quad   w_0(\pm \iota_0)=u_0(\pm \iota_0).
\end{equation}
Let $s=\big(\frac{\beta_n}{u_0(\iota_0)}\big)^{\frac1n} =\sqrt{1+\iota_0^2}$, then
 Lemma \ref{lm 5.1-1} with $x_0=\iota_0e$ implies that
\begin{equation}\label{a--0}
w_0(s\iota+\iota_0)=|\iota|^{-n} w_0(s\iota^{-1}+\iota_0) \quad {\rm and}\quad u_0(s\iota+\iota_0)=|\iota|^{-n} u_0(s\iota^{-1}+\iota_0).
\end{equation}
Note that
$$-\iota_0<s\iota^{-1}+\iota_0<\iota_0$$
is equivalent to that
 $$-\infty< s\iota+\iota_0< -\frac{1-\iota_0^2}{2\iota_0}. $$
Together with (\ref{a--}) and (\ref{a--0}), we obtain that
 \begin{equation}\label{a---}
w_0(\iota)>u_0(\iota) \quad {\rm for}\  \iota\in (-\infty, -\iota_1), \qquad w_0(\pm \iota_1)=u_0(\pm \iota_1),
\end{equation}
where
$$\iota_1=\frac{1-\iota_0^2}{2\iota_0}. $$

 If $\iota_0\geq 1$, then $0\geq \frac{1-\iota_0^2}{2\iota_0}>-\iota_0$ and  $w_0(\iota_1)>u_0(\iota_1)$ by (\ref{a--}),
 while a contradiction comes from (\ref{a---}), i.e.  $w_0(\iota_1)=u_0(\iota_1)$.

If  $\frac1{\sqrt{3}}<\iota_0<1$, then $\iota_1=\frac{1-\iota_0^2}{2t_0}\in(0,\iota_0)$,  $w_0(\iota_1)>u_0(\iota_1)$ by (\ref{a--}),
 while a contradiction comes from (\ref{a---}), i.e.  $w_0(\iota_1)=u_0(\iota_1)$.

 If $\iota_0= \frac1{\sqrt{3}}$, then $\iota_1=\iota_0$ and from  (\ref{a--}) and (\ref{a---}), we obtain that
 $$
w_0(\iota)>u_0(\iota) \quad {\rm for}\  \iota\in (-\infty, 0)\setminus\{\iota_0\}, \qquad w_0(\iota_0)=u_0(\iota_0)
$$
i.e.
$$
w (x)>u_1(x) \quad {\rm for}\  x\in\R^n\setminus  \{|z|=\iota_0\}, \qquad w_0(z)=u_0(z)\quad {\rm for}\ \, |z|=\iota_0,
$$
 which contradicts with (\ref{res-1}).

 If $0<\iota_0< \frac1{\sqrt{3}}$, then $ \iota_1=\frac{1-\iota_0^2}{2\iota_0}> \iota_0$ and  $w_0(-\iota_1)>u_0(-\iota_1)$ by (\ref{a--}),  a contradicts the fact that $w_0(-\iota_1)=u_0(-\iota_1)$ in \eqref{a---}.
 Now   we take $s=\sqrt{1+\iota_1^2}$, then Lemma \ref{lm 5.1-1} with $x_0=\iota_1 e$  gives the equalities
 $$w_0(s\iota+\iota_1)=\iota^{-n} w_0(s\iota^{-1}+\iota_1) \quad {\rm and}\quad u_0(s\iota+\iota_1)=|\iota|^{-n} u_0(s\iota^{-1}+\iota_1). $$
 Together with (\ref{a---}), we have that
  \begin{equation}\label{a--1}
 w_0(\iota)>u_0(\iota)\quad{\rm for}\ \, -\frac{1-\iota_1^2}{2\iota_1}< \iota<\iota_1,\qquad  w_0(\iota_1)=u_0(\iota_1),
 \end{equation}
  where $0<\iota_0<\iota_1$. However, this   contradicts  the equality $w_0(\iota_0)=u_0(\iota_0)$ in  (\ref{a--}).





When $w(0)<\beta_n$, it is similar to get a contradiction.    \hfill$\Box$\medskip

  \noindent{\bf Proof of Proposition \ref{teo 1-res}. } From Lemma \ref{lm 5.2-1}, we assume that $w$  is a solution of   \eqref{eq 1.1} verifying (\ref{res-1}) and
  $$ w(0)= w_\infty=\beta_n. $$
 We need to show that $w=u_1$ in $\R^n$.
We recall the notations
$$w_0(\iota)=w(\iota e)\quad{\rm and}\quad  u_0(\iota)=u_1(\iota e)\ \ {\rm for}\ \, \iota\in\R,$$
 where $e$ is any unit vector in $\R^n$.

 Assume that there exists $\iota_0>0$ such that
 \begin{equation}\label{b--0}
w_0(\iota_0)>u_0(\iota_0).
\end{equation}
  Since    $w(0)= w_\infty=\beta_n$, it follows by Lemma \ref{lm 5.1-1} with $s=1$, $t_0=0$ and $x=te$ that
 \begin{equation}\label{b--00}
 w_0(\iota)=|\iota |^{-n} v(\iota^{-1})\quad{\rm and}\quad  u_0(\iota)=|\iota|^{-n} u_0(\iota^{-1}),\quad \forall\, \iota\in\R.
 \end{equation}

Denote
 $$\tilde w_0(\iota)= |\iota|^{-n} w_0(\iota^{-1}+\iota_0) \quad{\rm and}\quad  \tilde  u_0(\iota)=|\iota|^{-n} u_0(\iota^{-1}+\iota_0) .  $$
 A direct computation shows that
   \begin{align*}
\tilde  u_0\big( \frac{-\iota_0}{1+\iota_0^2}\big)  &= \Big(\frac{ 1+\iota_0^2 }{2\iota_0}\Big)^n\, \frac{\beta_n}{\Big(1+(\frac{1+\iota_0^2}{-\iota_0}+\iota_0)^2\Big)^{\frac n2} }
\\[1mm]&= \beta_n  \Big(\frac{ 1+\iota_0^2 }{2\iota_0}\Big)^n\,  \Big(1+(\frac{1}{\iota_0} )^2\Big)^{-\frac n2}
 \\[1mm]&=\beta_n  (1+\iota_0^2)^{\frac n2}
  \\[1mm]& =(1+\iota_0^2)^n  \tilde  u_0(\iota_0)
\end{align*}
 and by (\ref{b--00}),
   \begin{align*}
\tilde  w_0\big( \frac{-\iota_0}{1+\iota_0^2}\big)  &=  \Big(\frac{ 1+\iota_0^2 }{2\iota_0}\Big)^n\, w_0\big(\frac{1+\iota_0^2}{-\iota_0}+\iota_0\big)
\\[1mm]&= \Big(\frac{ 1+\iota_0^2 }{2\iota_0}\Big)^n\, w_0\big(\frac{1}{\iota_0} \big)
 =   (1+\iota_0^2)^n  \tilde  w_0(\iota_0).
\end{align*}

 Furthermore, we have that
 \begin{align*}
\tilde  u_0(0)  &= \lim_{\iota\to0^+}  |\iota|^{-n} u_0(\iota^{-1}+\iota_0)\\[1mm]&=\beta_n  \lim_{\iota\to0^+} |\iota|^{-n}  \frac1{\Big(1+(\iota^{-1}+\iota_0)^2\Big)^{\frac n2} }
 = \beta_n  \lim_{\iota\to0^+}  \frac1{\Big(\iota^2+(1+\iota \iota_0)^2\Big)^{\frac n2} } =\beta_n
\end{align*}
 and
  \begin{align*}
\tilde  u_0\Big(\frac{-2\iota_0}{(1+\iota_0^2)}\Big)  &=  \Big(\frac{ 1+\iota_0^2 }{2\iota_0}\Big)^n\, \frac{\beta_n}{\Big(1+(\frac{1+\iota_0^2}{-2\iota_0}+\iota_0)^2\Big)^{\frac n2} }
\\[1mm]&= \beta_n  \Big(\frac{ 1+\iota_0^2 }{2\iota_0}\Big)^n\, \frac{1}{\Big(1+(\frac{\iota_0^2-1)}{2\iota_0} )^2\Big)^{\frac n2} }
\\[1mm]& =\beta_n \Big(\frac{ 1+\iota_0^2 }{2\iota_0}\Big)^n\, \frac{ (2\iota_0)^n }{\Big(\iota_0^2+\iota_0^2-2\iota_0^2+1 \Big)^{\frac n2} }=\beta_n.
\end{align*}
Note that  $-\frac{\iota_0}{(1+\iota_0)^2}$ is the center of $\tilde u_0$.  Let $\tau_0$ be the center of $\tilde w_0$. As (\ref{b--0}),
we have that
 \begin{equation}\label{b--1}
\tilde  w_0(\tau_0)\geq\tilde  w_0\big(-\frac{\iota_0}{(1+\iota_0)^2}\big)>\tilde  u_0\big(-\frac{\iota_0}{(1+\iota_0)^2}\big).
 \end{equation}
Furthermore, we see that
\begin{align*}
(\tilde  w_0)_\infty &= \lim_{\iota\to+\infty}  |\iota|^{-n} w_0(\iota^{-1}+\iota_0)
\\[1mm]& =\lim_{t\to+\infty} |\iota|^{-n} \big(|\iota| \big)^{n} v\big((\iota^{-1}+\iota_0)e\big)=v(\iota_0e)=w_0(\iota_0)
\end{align*}
and
\begin{align*}
(\tilde u_0)_\infty &= \lim_{\iota\to+\infty}  |\iota|^{-n} u_0(\iota^{-1}+\iota_0)
\\[1mm]& =\lim_{\iota\to+\infty}  u_1\big((\iota^{-1}+\iota_0)e\big)=u_1(\iota_0e)=u_0(\iota_0),
\end{align*}
then
 \begin{equation}\label{b--2}
(\tilde  w_0)_\infty > (\tilde u_0)_\infty.
 \end{equation}
  From (\ref{b--1}) and (\ref{b--2}), we have that
  $$\tilde w_0(\tau_0) (\tilde w_0)_\infty> \tilde u_0\big(\frac{-\iota_0}{(1+\iota_0)^2}\big) (\tilde u_0)_\infty=\beta_n^2.$$
 This contradicts Lemma \ref{lm 5.2-1}. Therefore, we get that
 $$w_0(\iota)\leq u_0(\iota)\quad{\rm for}\ \,  \iota\in\R, $$
 which says that
  $$w(x)\leq u_1(x)\quad{\rm for}\ \,  x\in\R^n $$
 which, together with the restriction (\ref{res-1}), we obtain that $w=u_1$ in $\R^n$.\smallskip

 If   there exists $t_0>0$ such that  $v(t_0)<u_1(t_0)$, a similar argument can lead to a contradiction. We complete the proof.   \hfill$\Box$\medskip

  \subsection{Full classification}

     \begin{proposition}\label{pro5-1low}
     Let $\Lambda_n$ be given in (\ref{Norm 1}),  $u\in L^2(\R^n)$ be a nontrivial classical solution of (\ref{eq 1.1}) and $u_t$ with $t>0$ be defined in (\ref{sol form-0}),  then
  $$\|u\|_{L^2(\R^n)} \geq \Lambda_n\ \big(= \|u_1\|_{L^2(\R^n)}\big). $$
    \end{proposition}

  \noindent{\bf Proof. }  Combining  (\ref{Norm 1}) and (\ref{stand sol con}), we see that
\begin{align*}
  \Lambda_n  =\beta_n \pi^{\frac n4} \sqrt{ \frac{\Gamma(\frac n2)}{\Gamma(n)}}
 =  (4\pi)^{\frac n4}  e^{\frac n2 \psi(\frac n2) } \sqrt{ \frac{\Gamma(\frac n2)}{\Gamma(n)}}.
  \end{align*}

  From \cite[Theorem 3]{Beckner}, that for $\|f\|_{L^2(\R^n)}=1$, there holds that
    \begin{align*}
    \frac n2 \int_{\R^n}(\ln|\xi|)\hat{f}(\xi)^2 d\xi \geq \int_{\R^n} (\ln |f|) f^2dx+B_n,
   \end{align*}
  where
$$B_n=\frac n2\Big(\psi(\frac n2) +\frac12 \ln \pi -\frac1n \ln\big( \frac{\Gamma(n)}{\Gamma(\frac n2)} \big)  \Big).$$
Moreover,  the equality holds for $u=u_1$.

Generally, for $f\not\equiv0$,
 \begin{align*}
   \int_{\R^n} (2\ln|\xi|)\hat{f}(\xi)^2 d\xi \geq   \frac 4n \int_{\R^n} \big(\ln |f|- \ln(\|f\|_{L^2(\R^n)}) \big) f^2dx+\frac 4nB_n\|f\|_{L^2(\R^n)}^2.
   \end{align*}
Particularly,
 \begin{align*}
   \int_{\R^n} (2\ln|\xi|)\hat{u_1}(\xi)^2 d\xi =   \frac 4n \int_{\R^n} \big(\ln |u_1|- \ln(\|u_1\|_{L^2(\R^n)}) \big) u_1^2dx+\frac 4nB_n\|u_1\|_{L^2(\R^n)}^2.
   \end{align*}
  Thus,
\begin{align}\label{dkdk}
\int_{\R^n}(\cL_\Delta u )u  dx  =   \frac 4n \int_{\R^n}(2\ln|\xi|)\hat{u }^2 d\xi
  \geq  \int_{\R^n} \big(\ln u - \ln(\|u \|_{L^2(\R^n)}) \big) u ^2dx+ \frac 4n B_n\|u \|_{L^2(\R^n)}^2
 \end{align}
 and
\begin{align}\label{dkdk}
\int_{\R^n}(\cL_\Delta u_1)u_1 dx  =   \frac 4n \int_{\R^n}(2\ln|\xi|)\hat{u_1}^2 d\xi
  =  \int_{\R^n} \big(\ln u_1- \ln(\|u_1\|_{L^2(\R^n)}) \big) u_1^2dx+ \frac 4n B_n\|u_1\|_{L^2(\R^n)}^2.
 \end{align}

  Multiply $u $ of (\ref{eq 1.1})   and integrate over $\R^n$,
then
  \begin{align*}
\int_{\R^n}(\cL_\Delta u )u  dx &= \frac 4n  \int_{\R^n} u ^2 \ln u  dx,
 \end{align*}
which together with (\ref{dkdk}), implies that
 \begin{align*}
\frac 4n \int_{\R^n} \big(\ln u- \ln(\|u\|_{L^2(\R^n)}) \big) u^2dx+\frac 4nB_n\|u\|_{L^2(\R^n)}^2 &\leq \frac 4n  \int_{\R^n} u^2 \ln u dx,
 \end{align*}
 and
 \begin{align*}
\frac 4n \int_{\R^n} \big(\ln u_1- \ln(\|u_1\|_{L^2(\R^n)}) \big) u^2dx+\frac 4nB_n\|u_1\|_{L^2(\R^n)}^2 = \frac 4n  \int_{\R^n} u_1^2 \ln u_1 dx,
 \end{align*}
Thus, we obtain that
  \begin{align*}
\Big(B_n-   \ln(\|u\|_{L^2(\R^n)})  \Big) \|u\|_{L^2(\R^n)}^2\leq 0
 \end{align*}
 and
 \begin{align*}
\Big(B_n-   \ln(\|u_1\|_{L^2(\R^n)})  \Big) \|u_1\|_{L^2(\R^n)}^2=0,
 \end{align*}
which say that
 \begin{align*}
  \ln(\|u\|_{L^2(\R^n)})\geq   B_n\quad{\rm and}\quad   \ln(\|u_1\|_{L^2(\R^n)})=   B_n,
 \end{align*}
   that is
 \begin{align*}
 \|u\|_{L^2(\R^n)} \geq   e^{B_n}=(4\pi)^{\frac{n}{4}}e^{\frac n2 \psi(\frac n2)} \sqrt{ \frac{\Gamma(\frac n2)}{\Gamma(n)} }= \Lambda_n
 \end{align*}
 and
 \begin{align*}
 \|u_1\|_{L^2(\R^n)}=  e^{B_n} = \Lambda_n.
 \end{align*}
 The proof is complete. \hfill$\Box$\medskip

     \noindent{\bf Proof of Theorem \ref{teo 1}. }
     Let $w$ be a nonnegative nontrivial solution of \eqref{eq 1.1} with maximum point at the origin  and denote
     $$w_0(\iota)=w(\iota e)\quad {\rm  for}\ \iota\in\R. $$

     From Theorem \ref{teo 1sym},     $w_0$ is a radial symmetric about $0$ and satisfies
\begin{equation}\label{fascinating asum-0}
 \lim_{|x|\to+\infty} w_0(r)r^n=w_\infty\in(0,+\infty),
 \end{equation}
 which implies that   $$ \|w\|_{L^2(\R^n)} \in(0,+\infty).$$
      Remark that  $w_0(\iota)=w_0(-\iota)$  by the radial symmetry.

  From Proposition \ref{pro5-1low}, we have that
  $$\|w\|_{L^2(\R^n)}\geq   \Lambda_n. $$
When $\|w\|_{L^2(\R^n)}=   \Lambda_n, $ we are done from Proposition \ref{teo 1-res}.
Next we can assume that
$$ \|w\|_{L^2(\R^n)} \in(\Lambda_n,+\infty).$$

       Recall that the standard solution  solution  of \eqref{eq 1.1}
      $$u_1(x)=\beta_n\frac{1}{(1+|x|^2)^{\frac n2}},\quad\forall\,  x\in\R^n  $$
      and
          $$u_0(\iota)=\beta_n\frac{1}{(1+\iota^2)^{\frac n2}},\quad\forall\,  \iota\in\R.   $$

Applying the  scaling properties in Theorem \ref{cr 3.0} $(ii)$,    we can assume that
there exists $\iota_0>0$ such that
\begin{equation}\label{fascinating asum-1}
w_0(\iota_0)=u_0(\iota_0),\qquad  w_0(\iota)>u_0(\iota)\ \ {\rm for}\ \,   |\iota|<\iota_0
\end{equation}
and,  together with (\ref{fascinating asum-0}),  there exits $0<r_1<r_2\leq +\infty$ such that for $i=1,2$
\begin{equation}\label{fascinating asum-2}
w_0(r_i)=u_0(r_i)\quad{\rm and}\quad   w_0(\iota)<u_0(\iota)\quad {\rm for}\ \, \iota \in(r_1,r_2).
      \end{equation}

Now we  argue as in the proof of Lemma \ref{lm 5.2-1} to obtain a contradiction,  that is
$$w_0\geq u_0\quad {\rm in}\ \,\R. $$

 Indeed, let $s=\big(\frac{\beta_n}{u_0(\iota_0)}\big)^{\frac1n} =\sqrt{1+\iota_0^2}$, then
 Lemma \ref{lm 5.1-1} with $x_0=\iota_0e$ implies that
\begin{equation}\label{a--0-n}
w_0(s\iota+\iota_0)=|\iota|^{-n} w_0(s\iota^{-1}+\iota_0) \quad {\rm and}\quad u_0(s\iota+\iota_0)=|\iota|^{-n} u_0(s\iota^{-1}+\iota_0).
\end{equation}
Note that
$$-\iota_0<s\iota^{-1}+\iota_0<\iota_0$$
which is equivalent to that
 $$-\infty< s\iota+\iota_0< -\frac{1-\iota_0^2}{2\iota_0}. $$
Together with (\ref{fascinating asum-1}) and (\ref{a--0-n}), we obtain that
 \begin{equation}\label{a----n}
w_0(\iota)>u_0(\iota) \quad {\rm for}\  \iota\in (-\infty, -\iota_1), \qquad w_0(\pm \iota_1)=u_0(\pm \iota_1),
\end{equation}
where
$$\iota_1=\frac{1-\iota_0^2}{2\iota_0}. $$

 If $\iota_0\geq 1$, then $0\geq \frac{1-\iota_0^2}{2\iota_0}>-\iota_0$ and  $w_0(\iota_1)>u_0(\iota_1)$ by (\ref{fascinating asum-1}),
 while a contradiction comes from (\ref{a----n}), i.e.  $w_0(\iota_1)=u_0(\iota_1)$.

If  $\frac1{\sqrt{3}}<\iota_0<1$, then $\iota_1=\frac{1-\iota_0^2}{2t_0}\in(0,\iota_0)$,  $w_0(\iota_1)>u_0(\iota_1)$ by (\ref{fascinating asum-1}),
again a contradiction comes from (\ref{a----n}), i.e.  $w_0(\iota_1)=u_0(\iota_1)$.

 If $\iota_0= \frac1{\sqrt{3}}$, then $\iota_1=\iota_0$ and from  (\ref{fascinating asum-1}) and (\ref{a----n}), we obtain that
 $$
w_0(\iota)>u_0(\iota) \quad {\rm for}\  \iota\in (-\infty, 0)\setminus\{\iota_0\}, \qquad w_0(\iota_0)=u_0(\iota_0)
$$
i.e.
$$
w (x)>u_1(x) \quad {\rm for}\  x\in\R^n\setminus  \{|z|=\iota_0\}, \qquad w_0(z)=u_0(z)\quad {\rm for}\ \, |z|=\iota_0,
$$
 which contradicts with (\ref{fascinating asum-2}).

 If $0<\iota_0< \frac1{\sqrt{3}}$, then $ \iota_1=\frac{1-\iota_0^2}{2\iota_0}> \iota_0$ and  $w_0(-\iota_1)>u_0(-\iota_1)$ by (\ref{fascinating asum-1}),  which contradicts the fact that $w_0(-\iota_1)=u_0(-\iota_1)$ in \eqref{a----n}.
 Now   we take $s=\sqrt{1+\iota_1^2}$, then Lemma \ref{lm 5.1-1} with $x_0=\iota_1 e$  gives the equalities
 $$w_0(s\iota+\iota_1)=\iota^{-n} w_0(s\iota^{-1}+\iota_1) \quad {\rm and}\quad u_0(s\iota+\iota_1)=|\iota|^{-n} u_0(s\iota^{-1}+\iota_1). $$
 Together with (\ref{a----n}), we have that
 $$
 w_0(\iota)>u_0(\iota)\quad{\rm for}\ \, -\frac{1-\iota_1^2}{2\iota_1}< \iota<\iota_1,\qquad  w_0(\iota_1)=u_0(\iota_1),
$$
  where $0<\iota_0<\iota_1$. However, this   contradicts  the fact that $w_0(\iota_0)=u_0(\iota_0)$ in  (\ref{a--}).

  Therefore, $\|w\|_{L^2(\R^n)} \in(\Lambda_n,+\infty)$ doesn't occur.
  \hfill$\Box$\medskip

\noindent{\bf Proof of  Theorem  \ref{teo-ful}. } It follows by  Theorem \ref{teo0} and Theorem \ref{teo 1} directly.   \hfill$\Box$

    \setcounter{equation}{0}
\section{Nonexistence by Pohozaev identity}
  \begin{proposition}\label{pr 6.1}
  Assume that
   $f$ is continuous in $\R$ such that
  $\displaystyle\lim_{t\to0}t^{-\gamma}f(t)=0$ for some $\gamma>0$
and
  $F(t)=\int_0^tf(\tau)d\tau$ for $t\in\R$.

  Let   $u\in L^1(\R^n)\cap \cD^0(\R^n)$ be a nonnegative function verifying
$uf(u), F(u)  \in L^1(\R^n)$ and
    \begin{equation}\label{eq 6.1}
  \cL_\Delta u=f(u) \quad {\rm in}\  \R^n,
  \end{equation}

 Then
   \begin{equation}\label{poh 6.1}
  2n\int_{\R^n}  F(u) dx -n\int_{\R^n} u f(u) dx +2\int_{\R^n} u^2 dx=0.
  \end{equation}

\end{proposition}

  \begin{lemma}\label{lm 6.1}
  For any function $\phi\in \cD(\R^n)$,
  \begin{equation}\label{eq 6.1-0}
  \cL_\Delta (x\cdot \nabla \phi) -(x\cdot \nabla)  \cL_\Delta\phi =2 \phi\quad {\rm in}\   \cD'(\R^n)
  \end{equation}
    \end{lemma}
  \noindent{\bf Proof. } From the basic properties of Fourier transformation, we have that
  $$\cF(\partial_j \phi)={\bf i\,} \xi_j \hat{\phi}\quad {\rm and}\quad \cF(x_j \phi)={\bf i\,} \partial_j \hat{\phi},  $$
  where ${\bf i}^2=-1$.
  Thus, we have that
 \begin{align*}
\cF(x\cdot \nabla \phi)  = \sum^n_{j=1} \cF( \partial_j (x_j \phi) -\phi) ={\bf i\,}  \sum^n_{j=1}  \xi_j \cF(x_j \phi) -n \hat{\phi}
=-\xi\cdot \nabla \phi -n \hat{\phi}
 \end{align*}
 and
  \begin{align*}
\cF\big(\cL_\Delta (x\cdot \nabla \phi)\big)  =-(2\ln|\xi|) \big(\xi\cdot \nabla \phi +n \hat{\phi}\big),
 \end{align*}
thus,
\begin{align*}
\cF\big((x\cdot \nabla)  \cL_\Delta\phi\big)  &= -\xi\cdot \nabla \cF\big(\cL_\Delta\phi\big)  -n  \cF\big(\cL_\Delta\phi\big)
\\[1mm]&=
-\xi\cdot \nabla  \big((2\ln |\xi|)\,  \hat{\phi}\big)  -n   (2\ln |\xi|)\,  \hat{\phi}
\\[1mm]&=  -2\hat{\phi} - (2\ln |\xi|)\,  \xi\cdot \nabla \hat{\phi}-n   (2\ln |\xi|)\,  \hat{\phi}
\\[1mm]&=-2\hat{\phi} +\cF\big(\cL_\Delta (x\cdot \nabla \phi)\big) ,
 \end{align*}
 that is,
 \begin{align*}
\cF\big(\cL_\Delta (x\cdot \nabla \phi)\big) =\cF\big((x\cdot \nabla)  \cL_\Delta\phi\big)  +2\hat{\phi}.
 \end{align*}
We complete the proof. \hfill$\Box$ \medskip

\noindent{\bf Proof of Proposition \ref{pr 6.1}. } Let $\rho_\epsilon(x)=\epsilon^{-n} \rho_0(\epsilon^{-1} x)$ for $x\in\R^n$, where
$\rho_0\in C^\infty(\R^n)$ is a nonnegative function having the support in $B_1(0)$ and
$\int_{\R^n}\rho_0(x)dx=1$. Then $\rho_\epsilon$ is a mollifier    and $u_\epsilon =\rho_\epsilon \ast u\in \cD(\R^n)$.
 In the sense of distributional sense, from (\ref{eq 6.1})
 $$\cL_\Delta u_\epsilon =\cL_\Delta (\rho_\epsilon \ast u)=(\cL_\Delta u)\ast \rho_\epsilon) =\rho_\epsilon \ast (\cL_\Delta u)=\rho_\epsilon \ast f(u)$$
 and since $x\cdot \nabla u_\epsilon \in \cD(\R^n)$, we have that
  \begin{align*}
  \int_{\R^n} \cL_\Delta u_\epsilon \big(x\cdot \nabla u_\epsilon\big) dx =\int_{\R^n} f(u)  \rho_\epsilon \ast \big(x\cdot \nabla u_\epsilon\big) dx=\int_{\R^n} f(u)  \big(x\cdot \nabla u_\epsilon\big) dx
  \end{align*}
  and by Lemma \ref{lm 6.1}
  \begin{align*}
  \int_{\R^n}f(u)  \big(x\cdot \nabla u_\epsilon\big) dx&= \int_{\R^n} \cL_\Delta u  \big(x\cdot \nabla u_\epsilon\big) dx
  \\[1mm] & =\int_{\R^n} u  \cL_\Delta \big(x\cdot \nabla u_\epsilon\big) dx
    \\[1mm] &=\int_{\R^n} u  \big((x\cdot \nabla)  \cL_\Delta u_\epsilon +2 u_\epsilon \big) dx
  \\[1mm]&=\int_{\R^n} u  \big((x\cdot \nabla)   \rho_\epsilon\ast f(u) +2 u_\epsilon \big) dx,
  \end{align*}
 then passing to the limit, we  obtain  that
   \begin{align}\label{eq 6-1---00}
  \int_{\R^n}f(u)  \big(x\cdot \nabla u \big) dx =\int_{\R^n} u  \big((x\cdot \nabla)    f(u) +2 u  \big) dx.
  \end{align}

 Note that
    \begin{align*}
   \int_{\R^n}f(u)  \big(x\cdot \nabla u \big) dx= \int_{\R^n}   \big(x\cdot \nabla F(u) \big) dx=-n\int_{\R^n} F(u)dx
     \end{align*}
     and
   \begin{align*}
   \int_{\R^n}u  \big(x\cdot \nabla f(u) \big) dx=- \int_{\R^n}f(u) {\rm div}(xu) dx=n\int_{\R^n} F(u)dx-n\int_{\R^n}u f(u)dx,
     \end{align*}
  then  combining  (\ref{eq 6-1---00}), we derive that
    \begin{align*}
 -n\int_{\R^n} F(u)dx=n\int_{\R^n} F(u)dx-n\int_{\R^n}u f(u)dx+2\int_{\R^n}u^2dx,
  \end{align*}
that is (\ref{poh 6.1}).  \hfill$\Box$ \medskip

\noindent  {\bf Proof of Theorem  \ref{teo 1.1-k}. } 
The prime function of $ t\ln t$ are $\frac12 t^2 \ln t-\frac14 t^2$ and we would apply (\ref{poh 6.1}) with $f(t)=t\ln t$ and $F(t)=\frac12 t^2 \ln t-\frac14 t^2$, by our assumption that $uf(u), F(u)\in L^1(\R^n)$,
 then we derive that
   \begin{equation}\label{poh 6.1-non}
 2\big(1- \frac n4 k\big)\int_{\R^n} u^2 dx=0.
  \end{equation}
  So when $k\not=\frac4n$, we have that $\int_{\R^n} u^2 dx=0$, which is impossible by our assumption that $u$ is nonzero.  \hfill$\Box$ \medskip

    \begin{corollary}\label{cr 6.1}
Let $p\in(0,1]$ and $u\in L^{2}(\R^n)\cap L^{p+1}(\R^n)\cap \cD^0(\R^n)$ be a solution of
\begin{equation}\label{eq 1.1-p}
\cL_\Delta u= |u|^{p-1}u \quad \   {\rm in}\ \,    \R^n.
\end{equation}
Then $u\equiv 0$.
\end{corollary}
  {\bf Proof. } From (\ref{poh 6.1}), we see that
$$
 n\big(\frac{ 2}{p+1}-1\big) \int_{\R^n}  |u|^{p+1} dx  +2\int_{\R^n} u^2 dx=0,
$$
which implies that $u=0$ in $\R^n$. \hfill$\Box$ \medskip

    \bigskip

\bigskip

\noindent {\bf  Conflicts of interest:} The authors declare that they have no conflicts of interest regarding this work.

\medskip\medskip

\noindent{\small {\bf Acknowledgements:} H. Chen is supported by  NNSF of China (Nos. 12071189, 12361043)
and Jiangxi Natural Science Foundation (No. 20232ACB201001).
F. Zhou was supported by Science and Technology Commission of Shanghai Municipality (No. 22DZ2229014)
and also NNSF of China (No. 12071189).}

\end{document}